\documentclass[reqno,11pt]{amsart}
\usepackage{amsmath,amssymb}
\usepackage{graphicx}
\usepackage[top=1in,bottom=1in,left=1in,right=1in]{geometry}
\newtheorem{thm}{Theorem}[section]

\newtheorem{lem}[thm]{Lemma}
\newtheorem{prop}[thm]{Proposition}
\newtheorem{defn}{Definition}[section]
\theoremstyle{definition}
\newtheorem{rem}{Remark}[section]
\numberwithin{equation}{section}

\usepackage{color}

\newcommand{\eps}{\varepsilon}

\def\be{\begin{equation}}
\def\ee{\end{equation}}
\def\rife#1{(\ref{#1})}
\def\vep{\varepsilon}
\def\de{\delta}
\def\vfi{\varphi}
\def\dint{\int\!\!\!\!\int}

\begin{document}
\title[Diffusive Hamilton-Jacobi equations]
      {The profile of boundary gradient blow-up \\ for the diffusive Hamilton-Jacobi equation}


\author[Porretta]{Alessio Porretta}%
\address{Universit\`a di Roma Tor Vergata,
Dipartimento di Matematica,
Via della Ricerca Scientifica 1,
00133 Roma, Italia}
\email{porretta@mat.uniroma2.it}

\author[Souplet]{Philippe Souplet}%
\address{Universit\'e Paris 13, Sorbonne Paris Cit\'e,
CNRS UMR 7539, Laboratoire Analyse, G\'{e}om\'{e}trie et Applications,
93430 Villetaneuse, France}
\email{souplet@math.univ-paris13.fr}

\begin{abstract}
We consider the diffusive Hamilton-Jacobi equation 
$$u_t-\Delta u=|\nabla u|^p,$$
with Dirichlet boundary conditions in two space dimensions,
which is a typical model-case in the theory of parabolic PDEs
and also arises in the KPZ model of growing interfaces. 
For $p>2$, solutions may develop gradient singularities on the boundary in finite time,
and examples of single-point gradient blowup on the boundary are known,
but the space-profile in the tangential direction has remained  a completely open problem.
In the parameter range $2<p\le 3$, for the case of a flat boundary and an isolated singularity at the origin,
we give an answer to this question, obtaining the precise final asymptotic profile, under the form
$$u_y(x,y,T) \sim d_p\Bigl[y+C|x|^{2(p-1)/(p-2)}\Bigr]^{-1/(p-1)},\quad\hbox{ as $(x,y)\to (0,0)$.}$$
Interestingly, this result  displays a new phenomenon of strong anisotropy of the profile, 
quite different from what is observed in other blowup problems for nonlinear parabolic equations,
with the exponents $1/(p-1)$ in the normal direction
$y$ and $2/(p-2)$ in the tangential direction $x$.
Furthermore, the tangential profile violates the (self-similar) scale invariance of the equation, 
whereas the normal profile remains self-similar.
\end{abstract}

\maketitle

\tableofcontents


\section{Introduction}

\subsection{Background}
This article is devoted to the qualitative study of solutions of the diffusive Hamilton-Jacobi equation
\begin{equation}\label{VHJ}
u_t-\Delta u=|\nabla u|^p.
\end{equation}
Beside being the viscosity approximation of Hamilton-Jacobi type
equations from stochastic control theory \cite{L82}, equation (\ref{VHJ}) is involved 
in certain physical models, for instance of ballistic deposition processes,
were it describes the evolution of the profile of a growing interface.
It is actually the deterministic version of the well-known
Kardar-Parisi-Zhang (KPZ) equation (see \cite{KPZ86} and cf. also Krug and Spohn \cite{KS88}).
In its stochastic version, it has undergone spectacular development recently with the work of M. Hairer~\cite{Hai13}.
Finally, equation (\ref{VHJ}) is a typical model-case in the theory of parabolic PDEs.
Indeed it is the simplest example of a parabolic equation with a nonlinearity depending on the first-order
spatial derivatives of $u$. As such, it is important to study its qualitative properties.

Equation (\ref{VHJ}) has been intensively studied in the past twenty years, and it is well known that two fundamentally different situations occur.
If the equation is considered in the whole space $\mathbb{R}^n$ (with, say, bounded $C^1$ initial conditions), then
all solutions exist globally in the classical sense and remain bounded in $C^1$; see e.g. \cite{AB98,BL99,BSW02,LS03,BKL04,GGK03,G05,S07,QS07},
where the large time behavior is investigated in detail).
At the opposite, if the equation is posed on a 
 domain with homogeneous Dirichlet boundary conditions, then for $p>2$ and suitably large initial data,
the local classical solution develops singularities in finite time.
These singularities are of gradient blowup type (GBU), the function~$u$~itself remaining bounded, and are located 
on some part of the boundary; see~e.g.~\cite{ABG89,FL94,CG96,A96,S02,ARS04,HM04,BaD04,SV06,SZ06,QS07,GH08,LS}
and the references therein.

For the classical blowup problem associated with the nonlinear heat equation
\begin{equation}\label{NLH}
u_t-\Delta u=u^p,
\end{equation}
a considerably developed theory is available for the description of the asymptotic profile of the solution
near a finite time singularity   (see \cite{QS07} and the references therein). In comparison, very little is known for equation (\ref{VHJ}).
In particular, in the case of an isolated boundary singularity, {\bf the final blowup profile of $\nabla u$
in the tangential direction is completely unknown.}\footnote{It is only known that 
 $|\nabla u(X,t)|\le C [{\rm dist}(X,\partial\Omega)]^{-1/(p-1)}$ (see \cite{ABG89, ARS04, SZ06, LS}), 
 which gives an upper estimate in the normal direction but provides no information on how the profile is damped away from the
 point of singularity along the boundary.}

\subsection{Main result: final gradient blowup profile near an isolated boundary singularity}

The goal of this paper is to fill this gap by giving a substantial contribution to this question.
In the range of exponents $2<p\le 3$, we will give a sharp description of 
the final blowup profile of $\nabla u$ near an isolated boundary singularity (in both normal and tangential directions).
Since the question is quite involved, we shall restrict ourselves to a rather simple setting, 
but which captures the essence of the problem. Namely, we consider the two-dimensional case,
where the domain is assumed to coincide locally with a half-plane near the point of singularity.
To this end, for given $\rho>0$, we set 
$$\omega_\rho=(-\rho,\rho)\times (0,\rho)\subset \mathbb{R}^2,\qquad  \omega_\rho^+=\omega_\rho\cap \{x>0\}.$$
Next we fix some $L,T>0$ and put
\begin{equation}\label{DefOmega1}
\omega=\omega_L,\quad \omega'=\omega_{L/2},
\end{equation}
\begin{equation}\label{DefOmega2}
Q_T:=\omega\times(0,T),\quad\Gamma_T= (-L,L)\times\{0\}\times (0,T).
\end{equation}

\begin{defn}
Let $L,T>0$ and let $u\in C^{2,1}(\overline\omega\times(0,T))$ be a nonnegative classical solution of (\ref{VHJ}) in $Q_T$,
with $u=0$ on $\Gamma_T$. We say that $u$ has an {\bf isolated gradient blowup point at~$(0,0,T)$} if
\begin{equation}\label{GBU0}
\limsup_{(x,y,t)\to (0,0,T)} |\nabla u(x,y,t)|=\infty
\end{equation}
and 
\begin{equation}\label{isolGBU}
\nabla u\ \hbox{ is bounded on $K\times (0,T)$ for any $K\subset\subset \overline\omega\setminus\{(0,0)\}$.}
\end{equation}
\end{defn}

If $u$ has an isolated gradient blowup point at $(0,0,T)$ then
we may define the {\bf final blowup profile} of $\nabla u$, given by
$$\nabla u(x,y,T):=\lim_{t\to T} \nabla u(x,y,t),
\quad\hbox{ for all $(x,y)\in \overline\omega'\setminus\{(0,0)\}.$}$$
Indeed the limit above  exists and is finite due to  (\ref{isolGBU}), as a consequence of standard parabolic estimates.
Our main result is the following.

\begin{thm}\label{thm:GBUprofile}
Assume
$$2<p\le 3.$$
Let $L,T>0$, let $u\in C^{1,2}(\overline\omega\times(0,T))$ be a nonnegative classical solution of (\ref{VHJ}) in $Q_T$,
with $u=0$ on $\Gamma_T$. Assume that $u$ has an isolated gradient blowup point at $(0,0,T)$ and that $u$ satisfies the monotonicity condition
\begin{equation}\label{MonotonicityHyp}
x{\hskip 0.4pt}u_x\le 0\quad\hbox{in $Q_T$.}
\end{equation}
Then there exist constants $C_1,C_2,C_3>0$, $\rho\in (0,L)$ (possibly depending on $u$) such that, for all $(x,y)\in ([-\rho,\rho]\times [0,\rho])\setminus\{(0,0)\}$,
the final blowup profile satisfies
\begin{equation}\label{mainEstimate}
d_p\Bigl[y+C_1 |x|^{2(p-1)/(p-2)}\Bigr]^{-\beta}-C_3 \le u_y(x,y,T) 
\le d_p\Bigl[y+C_2 |x|^{2(p-1)/(p-2)}\Bigr]^{-\beta} +C_3
\end{equation}
where
$$ \beta=1/(p-1)\quad\hbox{and}\quad d_p=\beta^\beta.$$
In particular, the final profile of the normal derivative on the boundary satisfies
$$C_4|x|^{-2/(p-2)}\le u_y(x,0,T)
\le C_5 |x|^{-2/(p-2)},$$
for all $0<|x|\le \rho$ and some $C_4,C_5>0$.
Also, for some $C>0$, we have
$$u\le C,\quad |u_x|\le C,\quad\hbox{ for all $(x,y)\in \omega'$}.$$

\end{thm}

   \medskip
 
 \hskip 1.5cm \includegraphics[width=13cm, height=7.5cm]{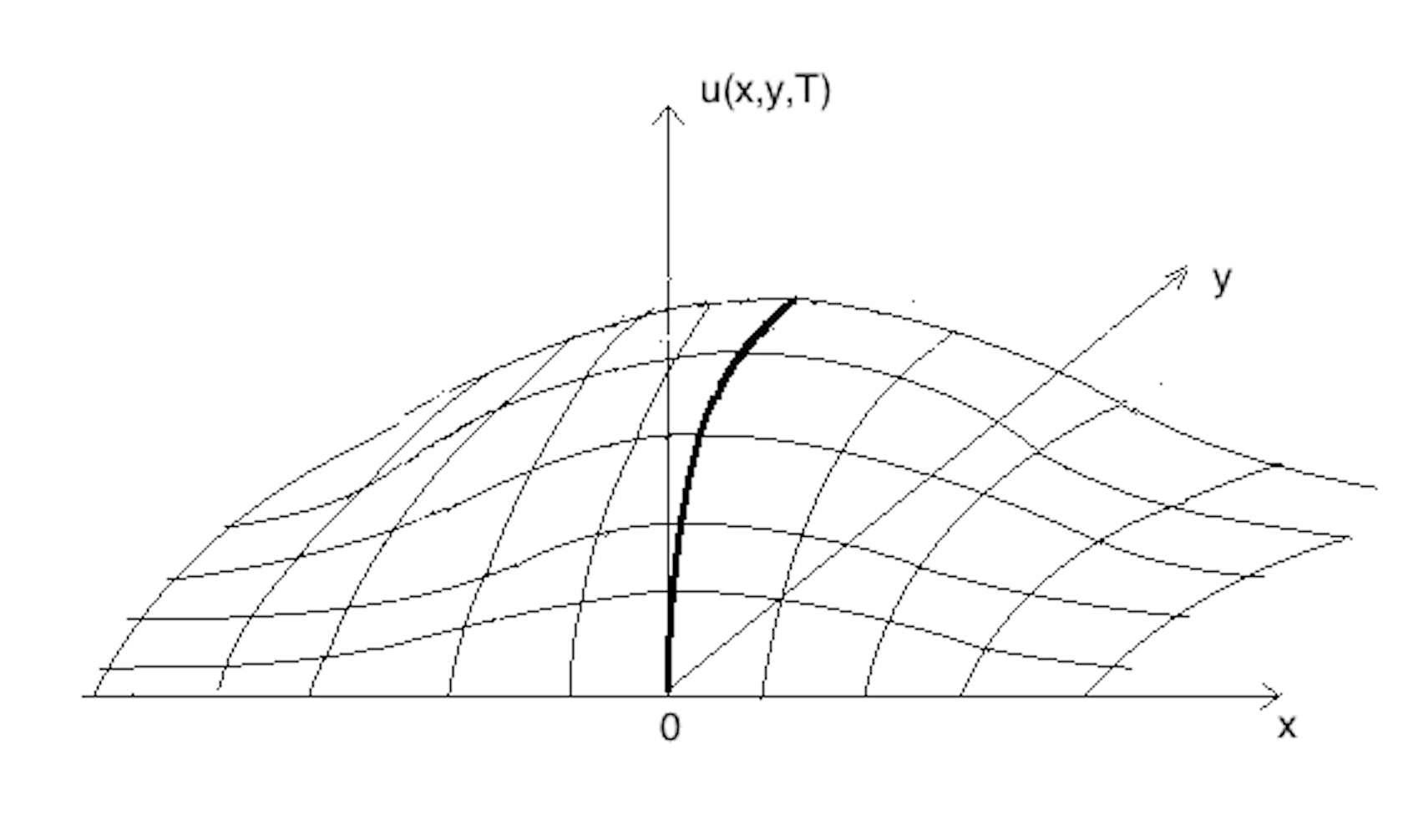}
 \vskip -2mm
  \centerline{\it Fig. 1: The shape of the final profile of $u$ near the origin.}
   \medskip
   
\subsection{Discussion and remarks}\hskip 3mm

\smallskip
(a) Interestingly, this result shows that the GBU profile is strongly {\bf anisotropic,}
i.e.~the exponents of the singularity profile in the normal and in the tangential directions are different,
respectively $1/(p-1)$ and $2/(p-2)$. 
Moreover, whereas the exponent of the normal profile obeys the natural scaling of the equation, 
the latter is violated by the tangential profile. Indeed, recall that equation (\ref{VHJ}) is invariant under the group of transformations
$$u \mapsto u_\lambda(x,y,t):=\lambda^m u(\lambda x,\lambda y, \lambda^2t)\quad\hbox{ with $m=(2-p)/(p-1)$, \quad for all $\lambda>0$},$$
whose gradient is given by $\nabla u_\lambda= \lambda^{1/(p-1)} \nabla u(\lambda x,\lambda y, \lambda^2t)$.

\smallskip
(b) As far as we know, no similar example of anisotropic, isolated blowup singularity is known in parabolic problems.
For the nonlinear heat equation~(\ref{NLH}), the stable blowup profile at an isolated blowup point is known to be
isotropic\footnote{this concerns blowup at an interior point -- actually only the whole space case is considered in these works; however no blowup can occur at a boundary point for equation (\ref{NLH}), at least in a convex domain}
 (see \cite{Ve93H, MZ97, MZ98, FMZ} and the references therein),
with 
$$u(X,T)\sim c(p)|X|^{-2/(p-1)}|\log |X||^{-1/(p-1)}\quad\hbox{as $X\to 0$}.$$
Here $X\in\mathbb{R}^n$ with $n\ge 2$ and $1<p<(n+2)/(n-2)$, and this profile occurs for instance for any symmetric, radially decreasing solution.

The case of the linear heat equation with nonlinear boundary conditions 
\begin{equation}\label{NLBC}
\begin{cases}
u_t-\Delta u=0\quad\hbox{in $\Omega\times(0,T) $}, \\
 \noalign{\vskip 1mm}
\displaystyle\frac{\partial u}{\partial\nu}=u^p\quad\hbox{on $\partial \Omega\times (0,T)$}\\
\end{cases}
\end{equation}
was studied in \cite{CF00, HY, Ha13, Ha15}.
Like for (\ref{VHJ}), this problem involves boundary singularities (however $u$ itself blows up). 
It was recently found in \cite{Ha13, Ha15} that for $\Omega=\mathbb{R}^2_+=\{(x,y);\ y>0\}$ under assumption (\ref{MonotonicityHyp}),
the singularity profile satisfies
$$u(x,y,T)\sim 
\begin{cases}
 y^{-1/(p-1)} &\quad\hbox{for $y\to 0$ with $|x|=O(y)$}\\
 \noalign{\vskip 1mm}
x^{-1/(p-1)}|\log x|^{-1/2(p-1)} &\quad\hbox{for $x\to 0$ and $y=0$.}
\end{cases}
$$
A similar result holds in dimension $n\ge 3$ if $1<p<n/(n-2)$. Note that this profile is only weakly anisotropic (by a logarithmic correction)
in comparison with (\ref{mainEstimate}).

 On the other hand we also observe that, unlike in problems (\ref{NLH}) and (\ref{NLBC}), 
the profile that we find for (\ref{VHJ}) is given by pure powers, without (e.g. logarithmic) corrections.
This situation seems typical of type II blow-up problems (see \cite{MM09} and cf. Remark~(c)).
\smallskip

(c) The exponent $2/(p-2)$ appears to be new in this problem. However, it is worth noting that, in some cases, the time rate of GBU involves
a related exponent.
Namely, for monotone in time solutions in 1 space dimension, we have \cite{GH08}:
\begin{equation}\label{TimeRate}
\|\nabla u(\cdot,t)\|_\infty\sim (T-t)^{-1/(p-2)}.
\end{equation}
However, the question of the time GBU rate is still open in 2 dimensions. 
Note that the rate (\ref{TimeRate}) corresponds to a type II blow-up, in the sense that this rate is more singular than what the natural scaling of the equation would suggest
(see \cite{GH08,QS07} for details).
A possible heuristic explanation of the appearance of the number $2/(p-2)$ in this problem, based on ideas of quasi-stationary approximation, is 
given in Section~6.

\smallskip
(d) It remains an open problem what is the actual tangential singularity exponent for $p>3$
-- see~Remark~\ref{remp3} for details.
Actually the lower estimate in (\ref{mainEstimate}) remains true for any $p>2$ (cf.~Theorem \ref{thm:GBUprofileLower}).
As for the upper estimate, for $p>3$, our method would allow to obtain
an estimate of the form in (\ref{mainEstimate}), with some power, greater than $2(p-1)/(p-2)$,
which could be explicitly computed in terms of $p$. However, due to the gap between the upper and
 lower estimates in this case, we are unable to determine the exponent of the actual profile.
Therefore, and in order not to further increase the technicality of the article,
we have refrained from expanding on this.
On the other hand, it might be possible to extend our results to more general (nonflat) domains
and to higher dimensions, at the expense of further complication.
But since the main goal of this work is to present a new phenomenon, we have decided to leave this aside.

\smallskip

(e) Actually, the upper estimate in (\ref{mainEstimate}) is satisfied by $u_y(x,y,t)$ for all $t<T$
 (this is a consequence of the proof, cf.~in particular formula (\ref{conclusionUpper})).

\smallskip

(f) By direct integration of (\ref{mainEstimate}) between $0$ and $y$, one easily obtains the corresponding estimate for the profile of the function $u$ itself
(whose shape is depicted in Fig.~1).
 In the course of the proof of Theorem \ref{thm:GBUprofile}, we also establish additional estimates, of possible independent interest. In particular, we show that
for any $p>2$, there holds 
$$|u_t|\le C,\qquad |u_x|\le C|x|,\qquad u_{xx}\ge -C$$
for $(x,y,t)$ close to $(0,0,T)$. Moreover, for $2<p\le 3$, we show that
$$-C\le u_{xx}(0,y,T) \le -c<0$$
for $y>0$ small (see Remark~\ref{rem51}). In particular, since $u_{xx}(x,0,T)=0$ for $x\ne 0$, we see that $u_{xx}(\cdot,T)$ is discontinuous near the origin.

\bigskip

\subsection{Existence of single-point gradient blow-up solutions}
In order to obtain solutions satisfying all the assumptions in Theorem \ref{thm:GBUprofile}
we now recall a result from  \cite{LS} concerning the initial-boundary value problem
\begin{eqnarray}
   &u_t-\Delta u=|\nabla u|^p,\ \ \ &x\in \Omega,\ t>0,\label{sys main 1}\\
   &u(x,t)=0,\ &x\in \partial \Omega,\ t>0,\label{sys main 2}\\
   &u(x,0)=u_0(x),\ &x\in \Omega.\label{sys main 3}
\end{eqnarray}
Here, it is assumed that
\begin{equation}\label{hypSPGBU}
\hbox{$\Omega\subset \mathbb{R}^2$ is a $C^{2+\alpha}$-smooth bounded domain,
$u_0\in C^1(\overline{\Omega})$ with $u_0\geq 0$ and ${u_0}_{|\partial \Omega}=0$.}
\end{equation}
 It follows from \cite[Theorem 10, p.~206]{F64} that problem (\ref{sys main 1})-(\ref{sys main 3}) admits a unique maximal,
nonnegative classical solution $u\in C^{2,1}(\overline{\Omega}\times(0,T))\cap
C^{1,0}(\overline{\Omega}\times[0,T))$, where $T=T(u_0)$ is the
maximal existence time.
Also, by the maximum principle, we immediately have
$$
    \|u(t)\|_\infty\leq \|u_0\|_\infty,\quad 0<t<T.
$$
On the other hand, by the Bernstein-type estimate in \cite{SZ06}, we know that
$$|\nabla u(t,X)|\le C[{\rm dist}(X,\partial\Omega)]^{-1/(p-1)},\quad \hbox{in $\Omega\times (0,T)$},$$
so that in particular GBU can take place only on $\partial\Omega$.
The following result was proved in \cite{LS}:
\medskip

{\bf Theorem A.} {\it
Assume (\ref{hypSPGBU}) and 
\begin{eqnarray}
    &&\hskip -1cm\hbox{$\Omega$ and $u_0$ are symmetric with respect to the line $x=0$}, \label{def dom 1}\\
    &&\hskip -1cm\hbox{$\Omega$ coincides locally near $0$ with the half-plane $\{y>0\}$ and is convex in the $x$-direction}, \label{def dom 2}\\
    &&\hskip -1cm\hbox{$xu_{0,x}\leq 0$ in $\Omega$}. \label{hyp ID 2}
\end{eqnarray}
If $u_0$ is suitably concentrated near the origin (see Remark~\ref{concentrated} below), then the solution of (\ref{sys main 1})-(\ref{sys main 3}) satisfies
$$\hbox{$T=T(u_0)<\infty$ and $\nabla u$ blows up only at the origin}$$
(i.e.  (\ref{GBU0}) is true and $\nabla u$ is bounded on $K\times (0,T)$ for any $K\subset\subset \overline\Omega\setminus\{(0,0)\}$).
}
\medskip

Also we note that, as a consequence of the assumptions of  Theorem~A,
we have $x{\hskip 0.4pt}u_x\leq 0$ in $\Omega\times (0,T)$.

\begin{rem}\label{concentrated}
As an example of data ``suitably concentrated near the origin'' for Theorem A, the following was given in  \cite{LS}:
\begin{equation}\label{def ID}
    u_0(x)=C\varepsilon^k\varphi\Bigr(\varepsilon^{-1}\sqrt{x^2+(y-\varepsilon)^2}\Bigl),
\end{equation}
where $k=(p-2)/(p-1)$, $C\geq C_0(p)>0$, 
$\varepsilon>0$ is sufficiently small, and
$\varphi\in C^\infty([0,\infty))$ satisfies
\begin{equation}\label{def ID2} 
    \varphi'\leq 0,\qquad \varphi(s)=1,\ s\leq 1/4, \qquad\varphi(s)=0,\
      s\geq 1/2.
\end{equation}
Note that these functions have support concentrated near the boundary point $(0,0)$
(and large derivatives for $\eps$ small).\end{rem}

\subsection{Ideas of proof}

The proof of Theorem \ref{thm:GBUprofile} is long and technical, and it requires to combine many ingredients.
Let us briefly describe the main ideas.

\smallskip
To establish the lower estimate, we start with an estimation of the normal derivative on the boundary,
which is obtained in three steps (see Fig.~2 in Section~3): we start from the vertical line $\{x=0\}$, where the 
precise final profile follows rather easily from ODE arguments. We then extend the lower estimate to the region 
above the curve
$y=K{\hskip 0.1pt}x^{2/(1-\beta)}$, by using a lower bound of $u_{xx}$ along horizontal segments.
The extension to the region below the curve $y=K{\hskip 0.1pt}x^{2/(1-\beta)}$ is then achieved
by means of a boundary Harnack-type estimate in suitable boxes connecting this curve to the boundary $\{y=0\}$.
Once $u_y$ is estimated from below on the boundary, the full lower estimate of $u_y$ is
obtained by suitable integration along vertical lines, plus some horizontal averaging
made possible by an estimate of the mixed derivative $u_{xy}$.

\smallskip
As for the proof of the upper estimate, it combines two ingredients. The first one is an auxiliary function of the form   
$$J(x,y,t)=u_x+kx\,(1+ y)\,y^{-(1-\beta)q}u^q,$$
with suitable parameters $k, q>0$,
which is shown to be nonpositive via the maximum principle.
The integration of the inequality $J\le 0$ along horizontal lines then yields a sharp upper estimate of H\"older type for $u$.
The second ingredient is a family of suitable regularizing barriers, which allow us to improve the H\"older estimate of $u$ to
a pointwise upper bound of $u_y$ on the boundary.
We note that rougher versions of both ingredients were used in \cite{LS},
in order to prove single-point GBU.{\footnote {\,The function $J$ in \cite{LS} was itself
motivated by a device from \cite{FML85}, where a function of the form $J(r,t)=r^{n-1}u_r+\eps r^n u^q$
was introduced to study the blowup of radial solutions of equation (\ref{NLH}).}
However, these ideas need to be considerably refined in order to obtain the sharp tangential GBU profile.
In particular, the derivation of the parabolic inequality satisfied by $J$ requires a delicate analysis in terms of
the auxiliary quantities 
$$\xi=y\frac{u_y}{u}\quad\hbox{ and }\quad\Theta=y(u_y)^{p-1}.$$
 This latter step turns out to require the restriction $p\le 3$ and leaves open the question whether the upper estimate remains true for $p>3$ as well (see~Remark~\ref{remp3}).

The rest of the paper is devoted to the proof of Theorem \ref{thm:GBUprofile}.
Some preliminary properties, mostly based on the maximum principle, are given in Section~2.
The lower estimate is established in Section~3. In Section~4, we construct the regularizing barriers.
In Section~5, the analysis of the parabolic inequality for the function $J$ is carried out, and the proof of the upper estimate is then completed.
Finally, a possible heuristic explanation of the appearance of the number $2/(p-2)$ in this problem, based on ideas of quasi-stationary approximation, is 
given in Section~6.

\section{Preliminary properties}

In the following propositions, we state a number of useful bounds and properties of the solution,
which will be used in the proof of the main result Theorem 1.1.
All the proofs will be given after the statements.
Here and in the rest of the paper, letters such as $C, C_1, c, \dots$ will denote possibly different positive constants,
whose dependence will be specified only when necessary.

We start with some simple bounds, which follow rather easily from the maximum principle.
Let us recall that $\omega,\omega', Q_T$ and $\Gamma_T$ are defined in \rife{DefOmega1}-\rife{DefOmega2}.

\begin{prop} \label{prop:maxprin}
Assume $p>2$, let $L,T>0$ and let $u\in C^{1,2}(\overline\omega\times(0,T))$ be a nonnegative classical solution of (\ref{VHJ}) in $Q_T$,
with $u=0$ on $\Gamma_T$. Assume that $u$ has an isolated gradient blowup point at~$(0,0,T)$.
\smallskip

(i) Then $u$ extends to a function 
\begin{equation}\label{extension}
u\in C^{1,2}(\tilde Q), \qquad\hbox{with $\tilde Q:=\bigl(\overline\omega'\times[T/2,T]\bigr)\setminus\{(0,0,T)\}$}.
\end{equation}
(This extension will still be denoted by $u$ without risk of confusion.)
\smallskip

(ii) There exists a constant $C>0$
(possibly depending on the solution $u$), such that $u$ satisfies the following bounds in $\tilde Q$:
\begin{eqnarray}
     &&|u_t|\le C \label{bound0}\\ 
     &&u_y\ge -C \label{bound2}\\
     &&u_{xx}\ge -C. \label{bound3}
\end{eqnarray}
 If, moreover, $u$ satisfies (\ref{MonotonicityHyp}), then we have
\begin{equation}\label{bound1}
 |u_x|\le C|x|
\end{equation}
in $\tilde Q$. 
\end{prop}

We next show that the gradient blowup does occur in a pointwise sense: $u_y$ becomes uniformly large near the blow-up time and the origin.

\begin{prop}\label{prop:blowup}
Assume $p>2$, let $L,T>0$ and let $u\in C^{1,2}(\overline\omega\times(0,T))$ be a nonnegative classical solution of (\ref{VHJ}) in $Q_T$,
with $u=0$ on $\Gamma_T$. Assume that $u$ has an isolated gradient blowup point at~$(0,0,T)$
and that $u$ satisfies the monotonicity condition (\ref{MonotonicityHyp}).
Then we have
\begin{equation}\label{limuy}
\lim_{{t\to T \atop (x,y)\to (0,0)}}u_y(x,y,t)=+\infty.
\end{equation}
\end{prop}

As a consequence of Proposition~\ref{prop:blowup}, there exists $0<\rho_0<\min(L/2,T/2)$ such that
\begin{equation}\label{uylarge}
u_y>1\quad\hbox{ in $\bigl(\overline\omega_{\rho_0}\times [T-{\rho_0},T]\bigl)\setminus\{(0,0,T)\}$.}
\end{equation}
 
We now give upper bounds which essentially follow by integrating in the vertical direction.

\begin{prop} \label{prop:upper1d}
Under the assumptions of Proposition \ref{prop:blowup},
there exist a constant $C>0$ and $\rho\in(0,\rho_0)$
(possibly depending on $u$), such that the solution $u$ satisfies
\begin{equation}\label{bound5a}
u_y(x,y,t)\le\bigl[(u_y)^{1-p}(x,0,t)+(p-1)y\bigr]^{-\beta}+2C y \quad\hbox{ in $\omega_\rho\times [T-\rho,T)$.}
\end{equation}
In particular, we have
\begin{eqnarray}
     &&u_y(x,y,t)\le d_p y^{-\beta}+2Cy \quad\hbox{ in $\omega_\rho\times [T-\rho,T)$} \label{bound5}\\
     &&u(x,y,t)\le c_p y^{1-\beta}+Cy^2  \quad\hbox{ in $\omega_\rho\times [T-\rho,T)$,}\label{bound6}
\end{eqnarray}
where 
\begin{equation}\label{defcpdp}
 \beta=\frac{1}{p-1}, \qquad d_p=\beta^\beta,\qquad c_p= (1-\beta)^{-1}d_p.
\end{equation}
\end{prop}

Our next result shows that similar lower bounds are true at $ x=0$
(of course they cannot be true for $x\neq 0$ in view of the profile eventually found in (\ref{mainEstimate})).

\begin{prop} \label{prop:lower1d}
Under the assumptions of Proposition \ref{prop:blowup},
there exist a constant $C_1>0$ and $\rho\in(0,\rho_0)$
(possibly depending on $u$) such that
\begin{equation}\label{bound7}
u_y(0,y,t)\ge \bigl[(u_y)^{1-p}(0,0,t)+(p-1)y\bigr]^{-\beta}-C_1.
\quad 0<y<\rho,\ \  T-\rho\le t<T.
\end{equation}
Moreover, we have
\begin{equation}\label{bound8}
u_y(0,y,T)\ge d_p y^{-\beta}-C_1, \quad 0<y<\rho
\end{equation}
and
\begin{equation}\label{bound9}
u(0,y,T)\ge c_p y^{1-\beta}-C_1y, \quad 0<y<\rho.
\end{equation}
\end{prop}

The following relationship between second order derivates, whose proof is rather delicate, 
will play an important role to establish the lower pointwise estimates in (\ref{mainEstimate}).

\begin{prop} \label{prop:uxy-uxx}
Under the assumptions of Proposition \ref{prop:blowup}, for any $\eta>0$, there exists a constant $C_\eta>0$
(possibly depending on $u$), such that the solution $u$ satisfies
\be\label{bound10}
u_{xy}\le \eta\, u_{xx}+C_\eta \qquad\hbox{in $\overline\omega_L^+\times [T/2,T]$.}
\ee
\end{prop}

\smallskip
We now turn to the proofs of the results that we have just stated.
\vskip2em

\noindent \textit{Proof of Proposition \ref{prop:maxprin}}.
$\bullet$ Property (\ref{extension}) is a consequence of standard parabolic estimates.

\smallskip
$\bullet$ Proof of (\ref{bound0})
and (\ref{bound2}). Set $z=u_t$ or $z=u_y$. Then $z\in C^{1,2}(\tilde Q)$ by parabolic regularity and it satisfies
\be\label{eqdrift}
z_t-\Delta z=A\cdot\nabla z,
\ee
where $A=p|\nabla u|^{p-2}\nabla u$.

Since $u_t=0$ on $\Gamma_T$, using (\ref{extension}), we see that the supremum of $|u_t|$ on the parabolic boundary of $\omega'\times[T/2,T)$ is finite.
Denoting this supremum by $C$, the maximum principle then guarantees
$$|u_t|\le C\quad\hbox{ in $\omega'\times[T/2,T)$,}$$
which implies (\ref{bound0}) in view of (\ref{extension}).

We can apply a similar reasoning to $u_y$. Since $u\ge 0$ and $u=0$ on $\Gamma_T$, we have $u_y\ge 0$ on $\Gamma_T$.
By (\ref{extension}), we see that the infimum of $u_y$ on the parabolic boundary of $\omega'\times[T/2,T)$ is finite.
Denoting this infimum by $-C$, the maximum principle then guarantees
$$u_y\ge -C\quad\hbox{ in $\omega'\times[T/2,T)$,}$$
which implies (\ref{bound2}) in view of (\ref{extension}).

\smallskip

$\bullet$ Proof of (\ref{bound3}). The function $Z:=u_{xx}\in C^{1,2}(\tilde Q)$ by parabolic regularity and it satisfies
\be\label{equxx}
\begin{split}
Z_t-\Delta Z
&= p\bigl[|\nabla u|^{p-2}\nabla u\cdot\nabla u_{x}\bigr]_{x} \\
&= A\cdot\nabla Z+p|\nabla u|^{p-2}|\nabla u_x|^2+p(p-2)|\nabla u|^{p-4}(\nabla u\cdot\nabla u_x)^2 \\
&\ge   A\cdot\nabla Z.
\end{split}
\ee
Since $Z=0$ on $\Gamma_T$, using (\ref{extension}) we see
that the infimum of $Z$ on the parabolic boundary of $\omega'\times[T/2,T)$ is finite.
It then follows from the maximum principle that $u_{xx}\ge -C$ in $\omega'\times[T/2,T)$,
which implies (\ref{bound3}) in view of (\ref{extension}).

\smallskip
$\bullet$ Proof of (\ref{bound1}). 
As a consequence of~(\ref{MonotonicityHyp}), we have
\begin{equation}\label{uxzero}
u_x(0,y,t)=0\qquad\hbox{for all $(y,t)\in \bigl([0,L/2]\times[T/2,T]\bigr)\setminus\{(0,T)\}$}.
\end{equation}Consequently, we get
$$u_x(x,y,t)=\int_0^x u_{xx}(t,s,y)\, ds\ge -Cx\qquad\hbox{in $\bigl(\omega'\times[T/2,T)\bigr)\cap\{x>0\}$,}$$
and a similar estimate for $x<0$.
This implies \rife{bound1}.
\hfill$\Box$

\vskip1em

In view of the proofs of Propositions \ref{prop:blowup}--\ref{prop:lower1d}, we prepare the following Lemma.

\begin{lem}\label{LemLower1}
Under the assumptions of Proposition \ref{prop:blowup}, we have
\begin{equation}\label{uyinfty}
\limsup_{t\to T}u_{y}(0,0,t)=+\infty
\end{equation}
and
\begin{equation}\label{bound8a}
u_y(0,y,T)\ge d_p y^{-\beta}-C_1, \quad 0<y<L/2.
\end{equation}
\end{lem}

\noindent \textit{Proof.}
As a consequence of (\ref{MonotonicityHyp}) and $u=0$ on $\Gamma_T$, we have 
\begin{equation}\label{uyxneg}
xu_{yx}(x,0,t)\le 0\quad\hbox{ for $0<|x|\le L/2$ and $t\in [T/2,T]$,}
\end{equation}
hence 
$$u_{y}(0,0,t)=\sup_{|x|\le L/2}u_{y}(x,0,t)\quad\hbox{ for $t\in [T/2,T)$.}$$
On the other hand, by (\ref{eqdrift}) and the maximum principle, for each $\tau\in (T/2,T)$, the maximum of $u_y$
 on $Q'_\tau:=\omega'\times (T/2,\tau)$ is attained on the parabolic boundary $\partial_PQ'_\tau$ of $Q'_\tau$.
 Moreover, by (\ref{isolGBU}), we have 
$$M_0=\sup_{(x,y,t)\in \Gamma'} u_y<\infty,
\quad\hbox{ where $\Gamma':= \partial_PQ'_T\setminus\bigl([-L/2,L/2]\times\{0\}\times [T/2,T)\bigr)$.}$$
Therefore,
 $$\sup_{Q'_\tau} u_y\le \max\Bigl(M_0,\sup_{t\in [T/2,\tau]} u_y(0,0,t)\Bigr).$$
By our assumption (\ref{GBU0}), the LHS goes to $\infty$ as $\tau\to T$
and (\ref{uyinfty}) follows.
\smallskip

Let us now prove (\ref{bound8a}).
By~(\ref{MonotonicityHyp}), we have
\begin{equation}\label{uxxzero}
u_{xx}(0,y,t)\le 0\qquad\hbox{for all $(y,t)\in (0,L/2]\times(T/2,T]$}.
\end{equation}
Also, we know from Proposition \ref{prop:maxprin}
 that $u_y\ge -C$ hence $|u_y|\le u_y+2C$ and that $|u_t|\le C$ in $\omega'\times (T/2,T)$.
Set $C_1:=2C+C^{1/p}$.
For $(y,t)\in (0,L/2)\times (T/2,T)$, using (\ref{uxzero}) and (\ref{uxxzero})  it follows that
$$ -u_{yy}(0,y,t)=[u_{xx}+|\nabla u|^p-u_t](0,y,t)\le |u_y(0,y,t)|^p+C\le (u_y(0,y,t)+C_1)^p.$$
Observing that $u_y(0,y,t)+C_1>0$ and integrating in $y$, we obtain
\begin{equation}\label{inequyyp}
(u_y(0,y,t)+C_1)^{1-p}\le (u_y(0,0,t)+C_1)^{1-p}+(p-1)y.
\end{equation}
By (\ref{extension}) and (\ref{uyinfty}), we deduce that
$$(u_y(0,y,T)+C_1)^{1-p}\le \liminf_{t\to T}(u_y(0,0,t)+C_1)^{1-p}+(p-1)y=(p-1)y\quad\hbox{ for $y\in(0,L/2)$},$$
which yields (\ref{bound8a}).
\hfill$\Box$
\vskip1em

\noindent \textit{Proof of Proposition \ref{prop:blowup}.}
Assume by contradiction that there exist a constant $K>0$ and a sequence  $(t_n, x_n, y_n)$ such that 
$$
(x_n, y_n, t_n) \to (0,0, T)\quad\hbox{ and }\quad u_y(x_n,y_n, t_n) \leq K\,.
$$
By (\ref{bound0}) and (\ref{bound3}), we have
\begin{equation}\label{bounduyyp}
u_{yy}+|u_y|^p\le  u_{yy}+|\nabla u|^p=u_t-u_{xx}\le 2C,
\end{equation}
hence in particular $u_{yy}\le 2C$.
Fix any $y\in (0,L/2)$. For $n$ large enough, we have $0<y_n<y$, hence
$$
u_y(x_n,y,t_n)\le u_y(x_n,y_n,t_n)+2C(y-y_n)\leq K + CL.
$$
Letting $n\to \infty$ and using (\ref{extension}), we get
$$
u_y(0,y,T)  \leq K + CL.
$$
This is in contradiction with \rife{bound8a}. 
\hfill$\Box$

\vskip1em
\noindent \textit{Proof of Proposition \ref{prop:upper1d}.}
Fix any $x\in (-L/2,L/2)$ and $t\in (T/2,T)$. By (\ref{uylarge}), there exists $\rho\in(0,\rho_0)$ such that
$$h(y) := u_y(x,y,t)-2Cy>0,\quad\hbox{ for $y\in (0,\rho)$}.$$
By (\ref{bounduyyp}), the function $h$ satisfies
$$h'+h^p =u_{yy}-2C+(u_y -2Cy)^p
\le u_{yy}-2C+(u_y)^p \le 0.$$
By integration, we obtain
$$h(y)\le\bigl[h^{1-p}(0)+(p-1)y\bigr]^{-\beta},\quad\hbox{ for $y\in (0,\rho)$},$$
hence (\ref{bound5a}) and in particular (\ref{bound5}). Property (\ref{bound6}) follows by further integration.
\hfill$\Box$

\vskip1em
\noindent \textit{Proof of Proposition \ref{prop:lower1d}}.
 Estimate (\ref{bound7}) is an immediate consequence of (\ref{inequyyp}).
As for (\ref{bound8}), it was already proved in Lemma~\ref{LemLower1}. Finally, (\ref{bound9}) follows from (\ref{bound8}) by integration. \hfill$\Box$

\vskip1em

We finally prove Proposition \ref{prop:uxy-uxx}.
\vskip1em

\noindent \textit{Proof of Proposition \ref{prop:uxy-uxx}}.  In view of estimate \rife{bound3}, there is no loss of generality if we only consider $\eta \leq 1$.  
First we recall from  \rife{equxx} that $u_{xx}$ satisfies
$$
(u_{xx})_t-\Delta u_{xx}
- A\cdot\nabla u_{xx}\geq 0\,,\qquad A=p |\nabla u|^{p-2}\nabla u \,.
$$
We compute the same equation for $u_{xy}$, and we get
\be\label{uxy1}
\begin{split}
(u_{xy})_t-\Delta u_{xy}
- A\cdot\nabla u_{xy} & = p|\nabla u|^{p-2}\Bigl[ \nabla u_y\cdot\nabla u_x+(p-2)\frac{(\nabla u\cdot\nabla u_x)(\nabla u\cdot\nabla u_y)}{|\nabla u|^2}\Bigr] 
\\
& = p|\nabla u|^{p-4}\Bigl[ |\nabla u|^2\nabla u_y\cdot\nabla u_x+(p-2) (\nabla u\cdot\nabla u_x)(\nabla u\cdot\nabla u_y)\Bigr].
\end{split}
\ee
Notice that this is justified close enough to the singularity, due to $u_y>0$ (cf.~(\ref{uyposQr}) below) and parabolic regularity.
Now, given $\eta \leq 1$, we consider the function
$$
z:= u_{xy}-\eta\, u_{xx}
$$
and the operator ${\mathcal L}(z)= z_t-\Delta z-A\cdot\nabla z$. On account of \rife{uxy1} we have
$$
{\mathcal L}(z)\leq F\qquad \hbox{in $ Q_T$}
$$
and
$$
F:= p|\nabla u|^{p-4}\Bigl[ |\nabla u|^2\nabla u_y\cdot\nabla u_x+(p-2) (\nabla u\cdot\nabla u_x)(\nabla u\cdot\nabla u_y)\Bigr]\,.
$$

We analyze now the sign of $F$ at large values of $z$. First of all, we develop each component of the scalar products and 
we find
\be\label{uxy2}
\begin{split}
\frac F{p|\nabla u|^{p-4}}  & = u_{xy} \left\{ u_{xx}( (p-1)(u_x)^2 + (u_y)^2) + u_{yy} ( (u_x)^2 + (p-1)(u_y)^2)\right\}
\\
&  \quad + (p-2) u_yu_xu_{yy}u_{xx} + (p-2) u_xu_yu_{xy}^2.
\end{split}
\ee
Due to Proposition \ref{prop:blowup}, we may choose $r>0$ small so that $u_y$ is sufficiently large in $Q'_r:=(0,r)^2\times (T-r,T)$. In particular, we may assume that
\be\label{uyposQr}
u_y >0\,,\qquad u_{yy}<0 \,,\qquad u_{xx}< - u_{yy}\qquad \hbox{in $Q'_r$,}
\ee
the last two inequalities coming from the equation $u_{xx}+ u_{yy} = -  |\nabla u|^p+ u_t $ together with the bounds \rife{bound0} and  \rife{bound3}.

In addition, since $Q'_r\subset \omega_L^+\times (0,T)$, we have $u_x\leq 0$ and therefore the last term in \rife{uxy2} is nonpositive. 
Dropping this term  we get
\begin{align*}
\frac F{p|\nabla u|^{p-4}}  & \leq  u_{xy} \left\{ u_{xx}( (p-1)(u_x)^2 + (u_y)^2) + u_{yy} ( (u_x)^2 + (p-1)(u_y)^2)\right\}
\\
&  \quad + (p-2) u_yu_xu_{yy}u_{xx} \,.
\end{align*}

Now,  since  $u_{xx}\geq -C$ (cf.~\rife{bound3}), at any point where $z\ge M_\eta:=\eta C$
we have $u_{xy} \geq M_\eta + \eta u_{xx}\geq 0$,
and since $u_{xx}\leq -u_{yy}$ we estimate
\be\label{uxy3}\begin{split}
\frac F{p|\nabla u|^{p-4}}   &  \leq  (p-2) \Bigl[ u_{xy}  u_{yy} ( (u_y)^2 - (u_x)^2) + u_yu_xu_{yy}u_{xx} \Bigr]
\\
& = (p-2) u_{yy} \Bigl[ u_{xy}   ( (u_y)^2 - (u_x)^2) + u_yu_x u_{xx} \Bigr]\,.
\end{split}
\ee
We conclude by noticing that the right hand side of \rife{uxy3} is negative if $u_y$ is large enough. Indeed, 
by Proposition \ref{prop:blowup} and (\ref{bound1}), we may chose $r$ so that
$$
u_y > \frac2\eta \|u_x\|_\infty \qquad \hbox{ in $Q'_r$.}
$$
Then, at any point of $Q'_r$ where $z\ge M_\eta$, we have
$$
u_{xy}   \geq  M_\eta  + \eta u_{xx}  
 = \eta (u_{xx}+  C) > \frac{2\|u_x\|_\infty}{u_y} (u_{xx}+  C)> -\frac{2u_x}{u_y} u_{xx}\,,
$$
hence
$$
u_{xy}   ( (u_y)^2 - (u_x)^2) + u_yu_x u_{xx} \geq  \frac12 u_{xy}  (u_y)^2 + u_yu_x u_{xx} 
=  u_y  \Bigl( \frac12 u_{xy}  u_y  +  u_x u_{xx}\Bigr) >0.
$$
So from \rife{uxy3} we get $F<0$ at any point of $Q'_r$ such that $z\ge \eta C$.
On the other hand, considering the parabolic boundary of $Q'_r$, we have $ z\leq \eta C$ at $\{x=0\}$ due to \rife{uxzero} and \rife{bound3}. 
At $\{ y=0\}$, we have $u_{xy}-\eta u_{xx}= u_{xy}\leq 0$ by (\ref{uyxneg}). On the rest of the lateral boundary, 
as well as at $t=T-r$, the function is bounded by some constant $C_\eta$. 
By the maximum principle applied to ${\mathcal L}$, we deduce that $ z\leq \max(\eta C,C_\eta)$ in $Q'_r$.
Therefore, the  bound \rife{bound10} is proved in $Q'_r$. In view of the regularity of $u$ outside of  $(0,0,T)$, the bound can of course be extended to 
$\omega_L^+ \times [T/2,T]$ up to an extra uniform constant.  
\hfill$\Box$

\section{Proof of main result: the lower estimate}
 
 In this section, we shall prove the following Theorem, which is valid for any $p>2$, and in particular implies the lower estimate in Theorem \ref{thm:GBUprofile}.

 \begin{thm}\label{thm:GBUprofileLower}
Assume $p>2$, let $L,T>0$ and recall the notation in (\ref{DefOmega1})-(\ref{DefOmega2}).
Let $u\in C^{1,2}(\overline\omega\times(0,T))$ be a nonnegative classical solution of (\ref{VHJ}) in $Q_T$,
with $u=0$ on $\Gamma_T$. Assume that $u$ has an isolated gradient blowup point at $(0,0,T)$ and that $u$ satisfies the monotonicity condition
\begin{equation}\label{MonotonicityHyp2}
x{\hskip 0.4pt}u_x\le 0\quad\hbox{in $Q_T$.}
\end{equation}
Then there exist constants $C_1,C_2>0$, $\rho\in (0,L)$ (possibly depending on $u$) such that, for all $(x,y)\in \overline\omega_\rho\setminus\{(0,0)\}$,
the final blowup profile satisfies
\begin{equation}\label{mainEstimateLower}
u_y(x,y,T)\ge d_p\Bigl[y+C_1 |x|^{2(p-1)/(p-2)}\Bigr]^{-\beta}-C_2,
\end{equation}
where
$$ \beta=1/(p-1)\quad\hbox{and}\quad d_p=\beta^\beta.$$
In particular, the final profile of the normal derivative on the boundary satisfies
\begin{equation}\label{mainEstimateLowerBoundary}
u_y(x,0,T) \ge C_3 |x|^{-2/(p-2)},
\end{equation}
for all $0<|x|\le \rho$ and some $C_3>0$.

\end{thm}

 In the rest of this section we denote  the final profile at the blow-up time by
 $$v:= u(\cdot,T)\in C^2(\overline\omega'\setminus\{(0,0)\}).$$

Theorem \ref{thm:GBUprofileLower} is proved in two steps. First, in Lemma~\ref{LemLower2}, we establish
the estimate of the normal derivative on the boundary 
(i.e. (\ref{mainEstimateLowerBoundary})).
To do so, the idea is as follows (see fig.~2 below): we start from the vertical line $\{x=0\}$, where the 
precise lower bound of the final profile $v$ is already known thanks to (\ref{bound9}).
We then extend the lower estimate of $v$ to the region $\Sigma_+$ above the curve
\begin{equation}\label{defParabolaSigma}
\Sigma_0=\bigl\{(x,y)\,: y=K{\hskip 0.1pt}x^{2/(1-\beta)}\bigr\},
\end{equation}
which plays an important role in our arguments.
This relies on the lower bound of $u_{xx}$ in Proposition~\ref{prop:maxprin}, used along horizontal segments.
This is not sufficient since the region $\Sigma_+$ does not touch the boundary $\{y=0\}$.
However, the extension to the region $\Sigma_-$ below the curve (\ref{defParabolaSigma}) can then be achieved
by using a Harnack-type estimate in suitable boxes connecting the  curve $\Sigma_0$ to the boundary $\{y=0\}$,
in terms of the distance to the boundary.

Finally, once the normal derivative is estimated on the boundary, the full lower estimate of $u_y$ is
obtained  (cf.~Lemma~\ref{LemLower3}) by suitable integration along vertical lines, plus some horizontal averaging
made possible by the estimate of the mixed derivative $u_{xy}$ given in Proposition \ref{prop:uxy-uxx}.

\begin{lem}\label{LemLower2}
 Under the assumptions of Theorem \ref{thm:GBUprofileLower},
there exist constants $c_0>0$ and $\rho\in (0,L)$ such that we have   
\be\label{vy0}
v_y (x,0) \geq c_0{\hskip 1pt} x^{- 2/(p-2)},\quad \hbox{ for $0<x<\rho$.}
\ee
\end{lem}

For the proof of Lemma \ref{LemLower2}, we shall use a well-known quantitative version of the Hopf Lemma (or 
boundary Harnack inequality) \cite{BC},
which we state in a suitably scale invariant form.

\begin{lem} \label{LemLowerBC}
Let $D_1$ be a $C^2$ domain such that 
$$ (-1, 1)\times (0,2) \subseteq D_1 \subseteq (-2, 2)\times (0,2).$$
For any $(x_0,y_0)\in  \mathbb{R}^2$ and $\lambda>0$, we set
$$D_\lambda:= (x_0,y_0) + \lambda D_1$$
and
$$\de_{D_\lambda}(x,y)={\rm dist}\bigl((x,y),\partial D_\lambda\bigr).$$
There exists $c_1>0$ depending only on $D_1$ such that for any $(x_0,y_0)\in  \mathbb{R}^2$,
any $\lambda>0$ and all $f\in L^\infty(D_\lambda)$, $f\ge 0$,
the solution $z$ of
\be\label{laplace}
\begin{cases}
-\Delta z=f & \hbox{in $D_\lambda$,}\\
z=0 & \hbox{on $\partial D_\lambda$}
\end{cases}
\ee
satisfies
$$
\frac{z(x,y)}{\de_{D_\lambda}(x,y)}\geq c_1\lambda^{-2} \dint_{D_\lambda} f(x',y')\de_{D_\lambda} (x',y') \,dx'dy',
\quad \hbox{ for all $(x,y)\in D_\lambda.$}
$$
\end{lem}

\vskip1em
\noindent \textit{Proof of Lemma \ref{LemLowerBC}.}
By translation invariance, we may assume $x_0=y_0=0$.
If $z$ solves (\ref{laplace}) in $D_\lambda$, then $Z(X,Y):=z(\lambda X,\lambda Y)$ solves
$$
\begin{cases}
-\Delta Z=f_\lambda(X,Y):=\lambda^2f(\lambda X,\lambda Y) & \hbox{in $D_1$,}\\
Z=0 & \hbox{on $\partial D_1$}.
\end{cases}
$$ 
The inequality for $\lambda =1$ is well known; see \cite{BC}. 
Using the fact that 
\be\label{scalingLambda}
\de_{D_1} (\lambda^{-1}x,\lambda^{-1}y)=\lambda^{-1}\de_{D_\lambda} (x,y)\quad \hbox{ for all $(x,y)\in D_\lambda,$}
\ee
and changing variables, it follows that
\begin{align*}
z(x,y)=Z(\lambda^{-1}x,\lambda^{-1}y)
& \geq c_1 \de_{D_1} (\lambda^{-1}x,\lambda^{-1}y)\dint_{D_1} f_\lambda(X',Y')\de_{D_1} (X',Y') \,dX'dY'\\
& = c_1 \de_{D_\lambda} (x,y)\dint_{D_1} f(\lambda X',\lambda Y')\de_{D_\lambda} (\lambda X',\lambda Y') \,dX'dY'\\
& = c_1 \lambda ^{-2}\de_{D_\lambda} (x,y)\dint_{D_\lambda} f(x',y')\de_{D_\lambda} (x',y') \,dx'dy',
\end{align*}
which proves the lemma.
 \hfill$\Box$

\vskip1em
\noindent \textit{Proof of Lemma \ref{LemLower2}.}
Starting from the lower estimate (\ref{bound9}) on $\{ x=0\}$, i.e.
\be\label{estiminfv00}
v(0,y) \geq c_p\, y^{1-\beta }- Cy\quad \hbox{ for $0<y<\rho$,}
\ee
the proof is done in three steps (cf. Fig.2 below).
\smallskip

{\bf Step 1.} {\it Lower estimate of $v$ in the region $\omega_\rho\cap\{y\ge K x^{2/(1-\beta)}\}$.} We claim that
there exist constants $K>0$ and $\rho\in (0,L)$ (depending on $v$) such that
\be\label{estiminfv0}
v(x,y) \geq \frac{c_p}{2}\,y^{1-\beta } \quad \hbox{ for $(x,y)\in \omega_\rho\cap\{y\ge K x^{2/(1-\beta)}\}$.}
\ee

Let $\rho$ be given by Proposition \ref{prop:upper1d}. 
Using the lower estimate (\ref{estiminfv00}) on $\{ x=0\}$, the fact that $v_x(0,y)=0$ and $v_{xx}\geq -C$ 
(cf. (\ref{uxzero}) and (\ref{bound3})) and Taylor's expansion, we obtain
$$
v(x,y) \geq c_p\, y^{1-\beta }- Cy- C \,x^2 \quad \hbox{ for $(x,y)\in \omega_\rho$.}
$$
hence
$$
v(x,y) \geq \bigl(c_p-CK^{\beta-1}-Cy^\beta\bigr)\, y^{1-\beta }  \quad \hbox{ for $(x,y)\in \omega_\rho\cap\{y> K x^{2/(1-\beta)}\}$.}
$$
The claim (\ref{estiminfv0}) follows by taking $\rho\le (c_p/4C)^{1/\beta}$ and $K\ge (c_p/4C)^{-1/(1-\beta)}$.
\smallskip

{\bf Step 2.} {\it Harnack-type estimate in suitable boxes near the boundary.}
We claim that there exist constants $c, \tilde c>0$ and $\rho\in (0,L)$ such that, for all $x\in (0,\rho/2)$ and all $\lambda\in (0,x/4)$,
\be\label{estiminfv1}
\frac{v(x,y)}{y}   \geq c\, \lambda^{1-2p}\, \left(\int_{x-\lambda}^x\int_0^{2\lambda} |v_y| \,dx'dy'\right)^{p}   - \tilde c \, \lambda
\quad \hbox{ for $0<y<\lambda$.}
\ee

By (\ref{bound3}) and (\ref{limuy}), reducing $\rho$ if necessary, we may assume that
$$-\Delta v \ge |\nabla v|^p-C\ge 0\quad \hbox{ in $\omega_\rho$.}$$
Let $D_1$ be a $C^2$ domain such that 
$$ (-1, 1)\times (0,2) \subseteq D_1 \subseteq (-2, 2)\times (0,2).$$
For given $x\in (0,\rho/2)$ and $\lambda\in (0,x/4)$, we set
$$
D=D_{x,\lambda}:= (x,0) + \lambda  D_1\subset (x/2,3x/2)\times (0,\rho)\subset  \omega_\rho\,.
$$
Observe that $-\Delta v \geq  f_{x,\lambda}:=|\nabla v|^p- C$ in $D$ with $f_{x,\lambda}\in L^\infty(D)$ and $f_{x,\lambda}\ge 0$.
Since $v\ge 0$, it follows from Lemma \ref{LemLowerBC} 
and the maximum principle that, for some constant $c_1>0$,
\be\label{estiminfv2}
\frac{v(x,y)}{\de_{D}(x,y)} \geq c_1 \lambda^{-2} \dint_{D} \left \{ |\nabla v|^p(x',y')- c \right\}\de_{D} (x',y')\,dx'dy'.
\ee

By H\"older's inequality, we have
$$
\dint_{D} |v_y| \,dx'dy'  \leq \left(\dint_{D} |v_y|^p\de_{D}(x',y')\,dx'dy'\right)^{\frac1p} 
\left(\dint_{D} \de_{D}^{-\frac1{p-1}}(x',y') \,dx'dy'\right)^{\frac{p-1}{p}}.
$$
Using (\ref{scalingLambda}), we see that
$$\dint_{D} \de_{D}^{-\frac1{p-1}}(x',y') \,dx'dy'
=\lambda^2\dint_{D_1} \de_{D}^{-\frac1{p-1}}(\lambda X',\lambda Y') \,dX'dY',
=\lambda^{2-\frac1{p-1}}\dint_{D_1} \de_{D_1}^{-\frac1{p-1}}(X',Y') \,dX'dY',$$
where the integral on the RHS is finite due to $1/(p-1)<1$ (see e.g. \cite{S02}).
Therefore, 
$$
\left(\dint_{D} |v_y| \,dx'dy'\right)^p \leq C \, \lambda^{2p-3} \dint_{D} |\nabla v|^p\de_{D}(x',y')\,  dx'dy'.
$$
Using also $\dint_{D} \de_{D}(x',y') \,dx'dy'=C\lambda^3$, we deduce
from (\ref{estiminfv2}) that
$$
\frac{v(x,y)}{\de_{D}(x,y)}   \geq c\, \lambda^{1-2p}\, \left(\dint_{D} |v_y| \,dx'dy'\right)^{p}   - \tilde c \, \lambda \,.
$$
Since $\de_{D}(x,y)= y$ for $0<y<\lambda$, the claim  (\ref{estiminfv1}) follows.
\smallskip

{\bf Step 3.} {\it Conclusion.}
Fix $x\in (0,\rho/2)$ (the case $x\in (-\rho/2,0)$ can be treated similarly).
We proceed to estimate from below the integral in (\ref{estiminfv1}).
To this end, we choose 
$$\lambda = K{\hskip 0.1pt}x^{2/(1-\beta)},$$
where $K$ is from (\ref{estiminfv0}).
Note that this implies $\lambda\in (0,x/4)$, taking a smaller $\rho$ if necessary.
By (\ref{estiminfv0}), we have
$$
\int_{x-\lambda}^x\int_0^{2\lambda} |v_y|\, dx'dy'\ge \int_{x-\lambda}^x\int_0^{2\lambda} v_y\, dx'dy'=  \int_{x-\lambda}^x v(x', 2\lambda) dx'
 \geq \frac{c_p}{2}\lambda (2\lambda)^{1-\beta}=c  \lambda^{2-\beta}.
$$
Combining this with (\ref{estiminfv1}) and  using $p\beta=\beta+1$   we obtain
$$
\frac{v(x,y)}{y} \geq c\, \lambda^{-\beta} -\tilde c\,\lambda 
\quad \hbox{ for $0<y<\lambda$,}
$$
which implies, reducing $\rho>0$ again if necessary,
\be\label{estiminfv3}
\frac{v(x,y)}{y}  \geq \frac{c}{2}\, x^{-2\beta/(1-\beta)}, \quad \hbox{ for $0<y<\lambda$.}
\ee
Since $2\beta/(1-\beta)= 2/(p-2)$,  letting $y\to 0$ we get \rife{vy0}.
 \hfill$\Box$
   \medskip
 
 \hskip 1.5cm \includegraphics[width=13cm, height=7.5cm]{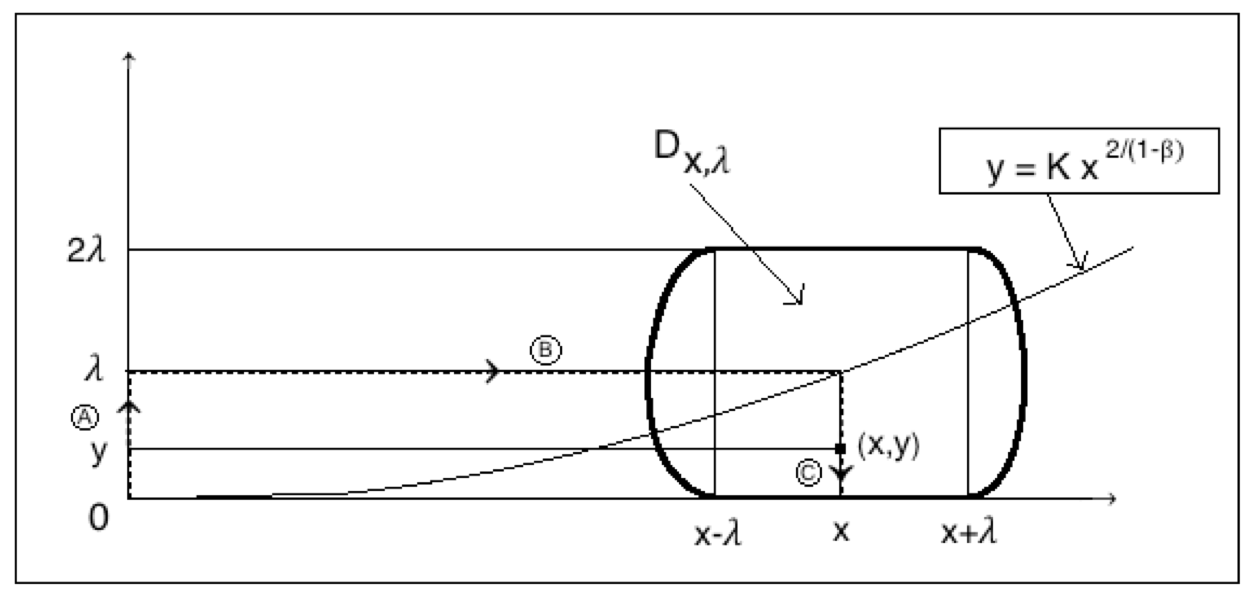}
 
  \centerline{\it Fig. 2: The scheme of the proof of Lemma 3.2.}
  \centerline{\it The marks A, B, C correspond to estimates 
  (\ref{estiminfv00}), (\ref{estiminfv0}) and (\ref{estiminfv3}), respectively}

 \medskip
The proof of Theorem \ref{thm:GBUprofileLower} is then completed by the following lemma. 
 \medskip

\begin{lem}\label{LemLower3}
 Under the assumptions of Theorem \ref{thm:GBUprofileLower},
for each $\vep\in(0,1)$, there exists constants $C>0$ and $\rho\in (0,L/2)$ such that 
\be\label{pre-media1}
v_y(x,y)   \geq \bigl[(v_y)^{1-p}((1+\vep)x, 0) +(p-1)y\big]^{-\beta} - C\quad\hbox{ in $\overline\omega_\rho\setminus\{(0,0)\}$.}
\ee
\end{lem}

\proof 
{\bf Step 1.} Let $\rho_0$ be given by (\ref{uylarge}) and take $\rho=\min(\rho_0,1)\le L/2$. We first claim that there exists a constant $A>0$ such that
\be\label{pre-media}
v_y(x,y) \ge \bigl[(v_y)^{1-p}(x,0)+(p-1)y\big]^{-\beta} - \int_0^y v_{xx}\,dy-A \quad\hbox{ in  $\overline\omega_\rho\setminus\{(0,0)\}$.}
\ee

From \rife{bound0}, \rife{bound1}, \rife{bound3} and (\ref{uylarge}), we know that
\be\label{lower1}
|u_t(\cdot,T)|\le C,\quad |v_x|\le C,\quad v_y>0, \quad v_{xx}\ge -C  \quad\hbox{ in  $\overline\omega_\rho\setminus\{(0,0)\}$}
\ee
for some constant  $C\ge 1$.
Let $A=3C$ and set 
$$
z(x,y):=v_y+\int_0^y v_{xx}\,dy+A.
$$
Observe that, for all $(x,y)\in \overline\omega_\rho\setminus\{(0,0)\}$
\be\label{lower2}
-z_y=-\Delta v \leq |\nabla v|^p + C.
\ee
Since $y\leq \rho\le 1$, using \rife{lower1}, we see that 
\be\label{lower3}
z = v_y+\Bigl(C+\int_0^y v_{xx}\,dy\Bigr)+2C\ge v_y+2C.
\ee
Therefore, by \rife{lower2}--\rife{lower3}, we obtain
$$
-z_y
\le [(v_y)^2+C^2]^{p/2}+C 
\le [v_y+C]^p+C^p 
\le [v_y+2C]^p 
\le z^p,
$$
hence
$$(z^{1-p})_y=-(p-1)z_yz^{-p}\le p-1\quad\hbox{ in $\overline\omega_\rho\setminus\{(0,0)\}$}.$$
After integration, it follows   that
$$z^{1-p}(x,y)\le z^{1-p}(x,0)+(p-1)y
\le (v_y)^{1-p}(x,0)+(p-1)y \quad\hbox{ in $\overline\omega_\rho\setminus\{(0,0)\}$,}
$$
hence the claim.

\smallskip

{\bf Step 2.} We may assume $0<x<\rho/2$ without loss of generality.
To prove \rife{pre-media1}, we now take the average of   inequality \rife{pre-media} and  we get
$$
\frac1{\vep x} \int_x^{(1+\vep)x}v_y(s,y)ds   \ge \frac1{\vep x} \int_x^{(1+\vep)x} \bigl[(v_y)^{1-p}(s,0)+(p-1)y\big]^{-\beta} ds- \frac1{\vep x} \int_x^{(1+\vep)x} \int_0^y v_{xx}\,dsdy-A\,.
$$
By \rife{uyxneg}, the function $x\mapsto v_y(x,0)$ is nonincreasing for $0<x<L$, so that we obtain
\be\label{lower4}
\begin{split}
\frac1{\vep x} \int_x^{(1+\vep)x}v_y(s,y)ds & 
\geq \bigl[(v_y)^{1-p}((1+\vep)x, 0) +(p-1)y\big]^{-\beta} - \int_0^y \frac1{\vep x} \int_x^{(1+\vep)x}  v_{xx}\,dsdy- A
\\ &
 \geq  \bigl[(v_y)^{1-p}((1+\vep)x, 0) +(p-1)y\big]^{-\beta} - C \frac y \vep- A,
 \end{split}
\ee
where we used \rife{bound1} in the last inequality.

Now we recall from Proposition \ref{prop:uxy-uxx} that, for some constant $C_1>0$ we have 
$$
(v_y-  v_x- C_1  x)_x \leq 0\qquad \hbox{in $\omega_L^+$.}
$$
Therefore, for some constant $C_2>0$, we get
\begin{align*}
v_y(x,y) & \geq \frac1{\vep x} \int_x^{(1+\vep)x}v_y(s,y)ds - \frac1{\vep x} \int_x^{(1+\vep)x}(  v_x(s,y)+ C_1  s)ds
+   v_x(x,y) + C_1 x
\\
  & \geq \frac1{\vep x} \int_x^{(1+\vep)x}v_y(s,y)ds - C_2 x\,,
\end{align*}
where again we also used \rife{bound1}. Jointly with inequality \rife{lower4}, we conclude that
$$
v_y(x,y)   \geq \bigl[(v_y)^{1-p}((1+\vep)x, 0) +(p-1)y\big]^{-\beta} - C_\vep
$$
i.e., \rife{pre-media1}.
 \hfill$\Box$
 
 \medskip
 
Finally, using Lemma \ref{LemLower3} with $\vep=1$, combined with Lemma \ref{LemLower2}, immediately yields \rife{mainEstimateLower}. 
The proof of Theorem \ref{thm:GBUprofileLower} is concluded.

\section{Regularizing barriers}

The following lemma shows that a   suitable local H\"older bound of exponent $1-\beta$,
near a boundary point,  actually guarantees a  bound for the normal derivative at this point.

\begin{lem}\label{barrier}
Let $p>2$, $r,d\in (0,1)$, $d<L$,   $t_0\in [0,T)$ and $x_0$ be such that $[x_0-r,x_0+r]\subset [-L/2,L/2]$.
Let 
$$D=(x_0-r,x_0+r)\times(0,d).
$$
There exist   constants $C_0=C_0(p)>0$ and $\eta_0=\eta_0(p,T)\in (0,1)$ 
 with the following property. Let $\eta\in(0,\eta_0)$ and $u\in C^{2,1}(\overline\omega\times(0,T))$ 
be a nonnegative classical solution of (\ref{VHJ}) in $Q_T$, with $u=0$ on $\Gamma_T$. If
\begin{equation}\label{hypcontroluc0}
    u(x,y,t)\leq c_py^{1-\beta}-  \kappa \frac{ y^2}2\ \ \ {\rm in}\ \
    \overline D \times[t_0,T),
\end{equation}
with $c_p=(1-\beta)^{-1}d_p$, 
$\kappa=C_0 \eta^{1-\beta}(r^2+ T-t_0)$ and
\be\label{toit}
u(x,d,t) \leq c_p \left[ (d+ \eta \,  (t-t_0)^{\frac1{1-\beta}}\,r^{\frac2{1-\beta}})^{1-\beta}- \eta^{1-\beta}\,  (t-t_0)\,r^2 \right] -  \kappa \frac{d^2}2  \ \ {\rm in}\ \
   [x_0-r,x_0+r]\times[t_0,T),
\, 
\ee
then 
\begin{equation}\label{conclcontroluc0}
u_y(x_0,0,t)\le  d_p \eta^{-\beta} \bigl(   (t-t_0)\,r^2 \bigr)^{-\frac\beta{1-\beta}}\ \ \ {\rm for }\ t\in ( t_0,T).
\end{equation} 

\end{lem}

 \begin{rem}
Lemma \ref{barrier} is an improvement with respect to \cite[Lemma 2.2]{LS}, where a similar result was proved,
but required a small constant instead of $c_p$ in assumption \rife{hypcontroluc0}.
This improvement is crucial in order to obtain the exact power in the upper estimate
of the GBU profile in the next section.
We note that assumption \rife{hypcontroluc0} is essentially sharp, since $u=c_py^{1-p}$ is an exact solution of (\ref{VHJ})
with $u=0$ and $u_y=\infty$ at $y=0$.
\end{rem}

\begin{proof}
Let us define the comparison function
\begin{eqnarray}\label{defzphi}
    z=z(x,y,t)=c_p\bigl[(y+\vfi(x,t))^{1-\beta}-\vfi^{1-\beta}(x,t)\bigr]-\kappa\frac{y^2}2\ \ \ {\rm in}\
    \overline D\times[t_0,T), 
\end{eqnarray}  
with
$$\vfi(x,t)=\eta\,   (t-t_0)^{\frac1{1-\beta}}\left(\frac{r^2-(x-x_0)^2}r\right)^{\frac2{1-\beta}},
$$
where $\eta>0$ and $\kappa\in (0,1)$. 
Let us denote by $C$ possibly different constants only depending on $p$ (often through the value of $\beta$). We first notice that there exists $C>0$ such that
\be\label{dervfi}
0\leq \vfi_t \leq C \eta^{1-\beta} r^2 \vfi^{\beta}\,,\qquad |\vfi_x|^2\leq C \eta^{ 1-\beta } (t-t_0) \vfi^{1+\beta}\,,\qquad |\vfi_{xx}|\leq C\eta^{1-\beta}   (t-t_0)\, \vfi^{\beta}\,.
\ee
Moreover,   if $\kappa$ is sufficiently small (depending only on $p,T$), we have
\begin{align}\label{kappasmall}
y(y+\vfi)^{\beta } \kappa  \, \le \kappa(1+T^{\frac{\beta}{1-\beta}})  < c_p(1-\beta) =d_p \,.
\end{align}
This implies 
$$
|z_y|^2 = d_p^2 (y+\vfi)^{-2\beta} \left( 1-\frac{\kappa y(y+\vfi)^{\beta}} d_p \right)^2 
\leq 
 d_p^2 (y+\vfi)^{-2\beta} 
$$
and since
$$
|z_x|^2 = d_p^2  \left(  (y+\vfi)^{-\beta}-\vfi^{-\beta}\right)^2 \vfi_x^2
\leq  d_p^2  \vfi^{-2\beta} \vfi_x^2,
$$
 we deduce
$$
|\nabla z|^p = (z_x^2+ z_y^2)^{\frac p2}\leq d_p^p\, (y+\vfi)^{-p\beta} \left[ 1+ (y+\vfi)^{2\beta}\vfi_x^2 \vfi^{-2\beta}\right]^{\frac p2}\,.
$$
Since, from \rife{dervfi}, we have $\vfi_x^2 \vfi^{-2\beta}\leq C\eta^{1-\beta} (t-t_0)\vfi^{1-\beta}$,  by taking $\eta_0=\eta_0(p,T)$ sufficiently small, it follows that $(y+\vfi)^{2\beta}\vfi_x^2 \vfi^{-2\beta}   \le C.$ Therefore we have
$$
|\nabla z|^p\leq  d_p^p\, (y+\vfi)^{-p\beta} \left[ 1+ C (y+\vfi)^{2\beta}\vfi_x^2 \vfi^{-2\beta}\right],
$$
which implies
\begin{align*}
|\nabla z|^p & \leq d_p^p\, (y+\vfi)^{-p\beta}  + C (y+\vfi)^{\beta-1}\vfi_x^2 \vfi^{-2\beta} 
\\ &
\leq  d_p^p\, (y+\vfi)^{-p\beta}  + C\vfi_x^2 \vfi^{-\beta-1} \,.
\end{align*}
Thus,  for $(x,t)\in D\times(t_0,T)$, we estimate
\begin{eqnarray*}
z_t-\Delta z- |\nabla z|^p  & \geq &  d_p\, [(y+\vfi)^{-\beta}-\vfi^{-\beta} ] (\vfi_t -\vfi_{xx}
)\\
 \noalign{\vskip 1mm} &&
+  d_p\, \beta[(y+\vfi)^{-\beta-1}-\vfi^{-\beta-1} ]  \vfi_x^2
\\
\noalign{\vskip 1mm} &&
  + d_p\,\beta(y+\vfi)^{-\beta-1} + \kappa
\\ \noalign{\vskip 1mm} &&
- d_p\, (y+\vfi)^{-p\beta}  - C\vfi_x^2 \vfi^{-\beta-1}.
\end{eqnarray*}
Using
$$
 d_p \beta(y+\vfi)^{-\beta-1} =  d_p^p\,  (y+\vfi)^{-p\beta}
$$
we deduce
$$
z_t-\Delta z- |\nabla z|^p    \geq   - C \vfi^{-\beta} (|\vfi_t| + |\vfi_{xx}|)
) 
-C \vfi^{-\beta-1}    \vfi_x^2 +\kappa 
$$
and thanks to  \rife{dervfi} we conclude
$$
z_t-\Delta z- |\nabla z|^p    \geq \kappa- C (r^2+  T-t_0) \eta^{1-\beta}\,.
$$
In particular, we have
\begin{eqnarray}\label{ineq:Pv negative}
    z_t-\Delta z\geq |\nabla z|^p\ \ \ {\rm in}\ D\times (t_0,T) 
\end{eqnarray}
with the choice $\kappa= C (r^2+  T-t_0) \eta^{1-\beta}$
  (which in turn guarantees (\ref{kappasmall}) for $\eta_0=\eta_0(p,T)$ small).

On $D\times \{  t_0\}$, as well as on
the lateral boundary part $\{x_0-r,x_0+r\}\times [0,d]\times (t_0,T)$, we have $\vfi=0$ and so 
$$
z= c_p y^{1-\beta} - \kappa  \frac{y^2}2 \geq u
$$
by \rife{hypcontroluc0}.  Next, on the part $[x_0-r,x_0+r]\times\{0\}\times (t_0,T)$, we 
have, for $t_0<t<T$,
$$
    u(\cdot,\cdot,t)=z(\cdot,\cdot,t)=0.
$$
On the remaining part $[x_0-r,x_0+r]\times\{d\}\times (t_0,T)$,  thanks to \rife{toit} 
  and to the fact that the expression in \rife{defzphi}
 is a decreasing function of $\varphi$, we have
$$
u(x,d,t)\leq c_p \left[ (d+ \eta \,  (t-t_0)^{\frac1{1-\beta}}\,r^{\frac2{1-\beta}})^{1-\beta}- \eta^{1-\beta}\,  (t-t_0)\,r^2 \right] - \kappa \frac{d^2}2 \leq z(x,d,t)\,.
$$
By the comparison principle, we deduce 
\begin{equation}\label{res:u leq v}
    u\leq z \ \ \ {\rm in}\
    \overline D\times[t_0,T). 
\end{equation}
In particular, we have
$$u_y(x_0,0,t)\le z_y(x_0,0,t)$$
  hence (\ref{conclcontroluc0}).
\end{proof}

\smallskip

\section{Proof of main result: the upper estimate}

We first establish the following

\begin{prop}\label{propJ} Let $p\in (2,3]$, $q>p-1$
 and let $u$ be as in Theorem 1.1. There exist  $k\in (0,1)$ and  $x_1, y_1, \sigma>0$ such that we have
$$
u_x+kxy^{-q(1-\beta)}(1+ y)u^q\leq 0 \ \ \ {\rm in}\ (0,x_1)\times(0,y_1)\times (  T-\sigma,T).
$$
\end{prop}

\begin{rem}\label{rem51}
It follows from Proposition~\ref{propJ} and estimate (\ref{bound9}) that there exists a constant $c>0$ such that, for $y>0$ small,
$$u_{xx}(0,y,T)=\lim_{x\to 0}\frac{u_x(x,y,T)}{x}\le -k(1+ y)\Bigl(\frac{u(0,y,T)}{y^{1-\beta}}\Bigr)^q \le -c<0.$$
\end{rem}

  \noindent \textit{Proof of Proposition \ref{propJ}.}
It is divided into several steps.

\smallskip

\textbf{Step 1. Preparations.} 
We consider the auxiliary function
\begin{equation*}
    J(x,y,t)=u_x+c(x)d(y)F(u)\ \ \ {\rm in}\ D\times(0,T),
\end{equation*}
where $D=(0,x_1)\times(0,y_1)$   and the smooth positive functions $c,d,F$ will be chosen below.
Our aim is to use the maximum principle to prove that,
for sufficiently small  $x_1, y_1, \sigma>0$, there holds
\begin{equation}\label{res:J negative}
    J\leq 0\ \ \ {\rm in}\ D\times[T  -\sigma,T).
\end{equation}
This will be done in the subsequent steps.

\smallskip

\textbf{Step 2. Derivation of a parabolic inequality for $J$.} 
The following basic computation was made in \cite{LS}.
For completeness and for the convenience of readers, we repeat it here.
We have
\begin{eqnarray*}
    J_t&=&u_{xt}+cdF'u_t,\\
    \Delta J&=&\Delta u_x+cdF'\Delta u+cdF''|\nabla u|^2
            +2c'dF'u_x+2cd'F'u_y+(cd''+ c''d)F. 
\end{eqnarray*}
Then we obtain
$$    J_t-\Delta J =(|\nabla u|^p)_x+cdF'|\nabla u|^p
            -2c'dF'u_x-2cd'F'u_y-cdF''|\nabla u|^2-(cd''+ c''d) F.$$
Using $u_x=J-cdF$, we write
\begin{eqnarray*}
    (|\nabla u|^p)_x&=&p|\nabla u|^{p-2}\nabla u\cdot\nabla u_x\\
            &=&p|\nabla u|^{p-2}\nabla u\cdot\nabla J
                    -p|\nabla u|^{p-2}\nabla u\cdot\nabla (cdF)\\
            &=&p|\nabla u|^{p-2}\nabla u\cdot\nabla J
                    -p|\nabla u|^{p-2}
                  (cdF'|\nabla u|^2+u_xc'dF+u_ycd'F),
\end{eqnarray*}
hence
\begin{eqnarray*}
    (|\nabla u|^p)_x&=&p|\nabla u|^{p-2}\nabla u\cdot\nabla J
                                   -pcdF'|\nabla u|^p \nonumber\\
            &&-pc'dF|\nabla u|^{p-2}J+pcc'd^2F^2|\nabla u|^{p-2}
                              -pcd'F|\nabla u|^{p-2}u_y. 
\end{eqnarray*}
We also have
\begin{equation*}
-2c'dF'u_x=-2c'dF'J+2cc'd^2FF'.
\end{equation*}
So we get
\begin{eqnarray}\label{eq:J}
    J_t-\Delta J&=&aJ+b\cdot\nabla J\nonumber\\
            &&-(p-1)cdF'|\nabla u|^p+pcc'd^2F^2|\nabla u|^{p-2}
                 -pcd'F|\nabla u|^{p-2}u_y
                              \nonumber\\
            &&+2cc'd^2FF'-2cd'F'u_y-cdF''|\nabla u|^2-(cd''+ c''d) F,
\end{eqnarray}
where
\begin{equation}\label{def a b}
    a=-pc'dF|\nabla u|^{p-2}-2c'dF' \ \ \ {\rm and}\ \ \
    b=p|\nabla u|^{p-2}\nabla u.
\end{equation}
Let
\begin{equation*}
    \mathcal{P}J=J_t-\Delta J-aJ-b\cdot\nabla J.
\end{equation*}
We can rewrite the above equality as follows
\begin{eqnarray}\label{ineq:J first estimate0}
    \frac{\mathcal{P}J}{cdF}& 
            &=-(p-1){F'\over F}|\nabla u|^p-{F''\over F}|\nabla u|^2  -{c''\over c}-{d''\over d} \\         
            &&\quad  -2{d'\over d}{F'\over F}u_y-p{d'\over d}|\nabla u|^{p-2}u_y
+pc'dF|\nabla u|^{p-2}+2c'dF'.\nonumber
\end{eqnarray}

\smallskip

\textbf{Step 3. Estimation of the RHS of (\ref{ineq:J first estimate0}).} 
We now specialize the previous computation to the following choice:
\begin{eqnarray}
    &&F(u)=u^q,\ \ \ \label{def F}\\
    &&d(y)=y^{-\gamma}\vfi(y),\label{def d}\\
    &&c(x)=kx,\label{def c}
\end{eqnarray}
where $k\in (0,1)$, $q>1$, $\gamma>0$ and $\varphi>0$ is a smooth function with $\varphi'\ge 0$.
Using 
$$
d'(y)=-\gamma y^{-\gamma-1}\vfi(y)+y^{-\gamma}\varphi'(y),
$$
$$
d''(y)=\gamma(\gamma+1) y^{-\gamma-2}\vfi(y)-2\gamma y^{-\gamma-1}\varphi'(y)+y^{-\gamma}\varphi''(y),
$$
the equality \rife{ineq:J first estimate0} implies 
\begin{eqnarray}\label{ineq:J first estimate}
    \frac{y^2\mathcal{P}J}{cdF}&=
            &-(p-1)q\frac{y^2|\nabla u|^p}u-q(q-1)\frac{y^2|\nabla u|^2}{u^2}-\gamma(\gamma+1) +2\gamma \frac{y\varphi'}{\varphi}-\frac{y^2\varphi''}{\varphi}\\         
            &&+2q\Bigl(\gamma-\frac{y\varphi'}{\varphi}\Bigr)\frac{yu_y}{u}+p\Bigl(\gamma-\frac{y\varphi'}{\varphi}\Bigr)y|\nabla u|^{p-2}u_y
+pky^{2-\gamma}\varphi u^q|\nabla u|^{p-2}+2kqu^{q-1}y^{2-\gamma}\varphi.\nonumber
\end{eqnarray}
Also, taking $\sigma,x_1,y_1\in(0,1)$ sufficiently small and setting   
$Q=(0,x_1)\times (0,y_1)\times (T-\sigma,T)$, we have,  by \rife{bound1}, \rife{uylarge}, \rife{bound5}, \rife {bound6},
$$
|u_x|\le Cx \qquad\hbox{in $Q$},  
$$
as well as
\be\label{u1d}
y\le u\le c_py^{1-\beta}+Cy^2,\quad 1\le  u_y\le d_py^{-\beta}+Cy\qquad\hbox{in $Q$}.
\ee
In particular we have, close enough to the singularity,
\begin{align*}
y|\nabla u|^{p-2}u_y
&=y(u_y)^{p-1}\Bigl[1+\Bigl(\frac{u_x}{u_y}\Bigr)^2\Bigr]^{(p-2)/2}\le y(u_y)^{p-1}
\Bigl[1+\frac{p-2}2\Bigl(\frac{u_x}{u_y}\Bigr)^2\Bigr]\\
&=y(u_y)^{p-1}+ \frac{p-2}2 y(u_y)^{p-3}(u_x)^2
\\
& \le y(u_y)^{p-1}+Cy^m\, x^2 
\end{align*}
with $m=\min(1,2\beta)$. In particular, for $p\leq 3$ (i.e. $\beta\geq \frac12$), we have
$$
 y|\nabla u|^{p-2}u_y\leq y(u_y)^{p-1}+ Cy\, x^2.
 $$
Similarly, using \rife{u1d} we estimate 
$$
pky^{2-\gamma}\varphi u^q|\nabla u|^{p-2}+2kqu^{q-1}y^{2-\gamma}\varphi\le kqC\vfi\,y^{(q-1)(1-\beta)+2-\gamma}
\leq q C\vfi\, y^{1+\beta  +q(1-\beta)-\gamma} 
$$
for   any $k\in (0,1)$.
Consequently, we get  from \rife{ineq:J first estimate}
\be\label{pregam}
\begin{split}
{y^2{\mathcal P}J\over cdF}
&\le-(p-1)q\frac{y^2(u_y)^p}{u} -q(q-1)\frac{y^2(u_y)^2}{u^2}-\gamma(\gamma+1)\\
&\quad  +2q\Bigl(\gamma-\frac{y\varphi'}{\varphi}\Bigr)\frac{yu_y}{u}+p\Bigl(\gamma-\frac{y\varphi'}{\varphi}\Bigr)y(u_y)^{p-1}\\
&\quad +2\gamma \frac{y\varphi'}{\varphi}-\frac{y^2\varphi''}{\varphi} 
+ Cy(x^2+ q\varphi y^{\beta +q(1-\beta)-\gamma})  
\qquad\hbox{ in $Q$}.
\end{split}
\ee

We will conclude Step 3 through the following lemma.

\begin{lem}\label{pjneg} 
Let $p\in (2,3]$, $q>p-1$, and take  $\gamma=q (1-\beta)$ and $\vfi(y)=1+ y$ in  \rife{def d}.  
There exist $x_1, y_1, \sigma>0$  sufficiently small such that,  
for any $k\in (0,1)$, we have  
$${\mathcal P}J\le 0 \quad\hbox{in $Q=(0,x_1)\times(0,y_1)\times (T-\sigma,T)$.}
$$
\end{lem}

  \noindent \textit{Proof of Lemma \ref{pjneg}.} To shorten notations, we set
$$
\xi=\xi(x,y,t)=\frac{yu_y}{u}\, \ge 0,\qquad \Theta=\Theta(x,y,t)= y(u_y)^{p-1}\, \ge 0, \qquad \psi=\psi(y)=\frac{y\varphi'}{\varphi} .
$$
Notice that $\psi(y)= \frac{ y}{1+ y}$ is small provided  $y_1$ is sufficiently small.

We wish to show that the right-hand side of \rife{pregam} is nonpositive in $Q$. To this purpose we distinguish two cases according to whether $\xi\leq 1-\beta$ or not.

\smallskip
\noindent{\hskip 4mm}{\bf Case 1: $\xi\leq 1-\beta$.}
\smallskip

By Young's inequality, we have
\be\label{add-rev1}
p\Bigl(\gamma-\psi(y)\Bigr)y(u_y)^{p-1} \leq  (p-1) q \frac{y^2(u_y)^p}{u}+ 
\Bigl(\gamma-\psi(y)\Bigr)^{p} \frac{u^{p-1}}{q^{p-1}y^{p-2}}\,.
\ee
Using \rife{u1d}, and $1-\beta= (p-2)/(p-1)$, we have
\begin{align*}   
\Bigl(\gamma-\psi(y)\Bigr)^{p} \frac{u^{p-1}}{q^{p-1}y^{p-2}} & \leq \gamma^p\Bigl(1-\frac{\psi(y)}\gamma\Bigr)^{p} \frac{(c_p y^{1-\beta}+ C y^2)^{p-1}}{q^{p-1}y^{p-2}}
\\
& \quad \leq  \gamma \Bigl(1-\frac{\psi(y)}\gamma\Bigr)^{p} \left( \frac{c_p\gamma}q\right)^{p-1}  \left(1+ C y^{1+\beta}\right)^{p-1}.
\end{align*}
The precise value of $c_p$ in \rife{defcpdp} and the choice $\gamma=q (1-\beta)$ then imply
$$
\Bigl(\gamma-\psi(y)\Bigr)^{p} \frac{u^{p-1}}{q^{p-1}y^{p-2}}   \leq 
\gamma\beta \Bigl(1-\frac{\psi(y)}\gamma\Bigr)^{p}  \left(1+ C y^{1+\beta}\right)^{p-1} \leq  q(1-\beta)\beta- p \beta \psi(y)  + C y^{1+\beta}
$$
for some $C$ (possibly depending on $q$). 
 Hence from \rife{add-rev1} we obtain
$$
p\Bigl(\gamma-\psi(y)\Bigr)y(u_y)^{p-1} \leq  (p-1) q \frac{y^2(u_y)^p}{u}+ 
q(1-\beta)\beta- p \beta \psi(y)  + C y^{1+\beta}\,.
$$ 
Therefore, using $\vfi''=0$, \rife{pregam} implies
\begin{align*}
{y^2{\mathcal P}J\over cdF}
&\le  -q(q-1)\xi^2 +2q\Bigl(\gamma-\psi(y)\Bigr)\xi \\
&\quad  - \gamma(\gamma+1)+ q(1-\beta)\beta- p \beta \psi(y) +  2\gamma \psi(y)+ Cy(x^2+ y^\beta).
\end{align*}
Now we remark that the function 
$$
\xi\mapsto  -q(q-1)\xi^2  +2q\Bigl(\gamma-\psi(y)\Bigr)\xi
$$
is increasing for $\xi \leq 1-\beta$  and $y$ sufficiently small, so we get
\begin{align*}
{y^2{\mathcal P}J\over cdF}
&\le  -q(q-1)(1-\beta)^2+2q\Bigl(\gamma-\psi(y)\Bigr)(1-\beta)   \\
&\quad  - \gamma(\gamma+1)+ q(1-\beta)\beta- p \beta \psi(y) +  2\gamma \psi(y)+ Cy(x^2+  y^\beta)
\end{align*}
which implies, using $\gamma=(1-\beta)q$,
$$
{y^2{\mathcal P}J\over cdF} \leq - p \beta \psi(y) + Cy(x^2+  y^\beta)=- p \beta \frac{ y}{1+ y} + Cy(x^2+ y^\beta)\,.
$$
Therefore, we have ${\mathcal P}J\leq 0$  provided $x_1$, $y_1$ are taken sufficiently small. 
\medskip

\noindent{\hskip 4mm}{\bf Case 2: $\xi > 1-\beta$.}
\smallskip

With the above notations, and using $\vfi''=0$,   \rife{pregam} can be written as
\be\label{pregam2}
\begin{split}
{y^2{\mathcal P}J\over cdF}
&\le-(p-1)q\xi\,\Theta -q(q-1)\xi^2-\gamma(\gamma+1) \\
&\quad  +2q\Bigl(\gamma-\psi(y)\Bigr)\xi+ p\bigl(\gamma-\psi(y)\bigr) \Theta \\
&\quad +2\gamma \psi(y)+C y(x^2+ y^\beta)
\end{split}
\ee
hence
$$
{y^2{\mathcal P}J\over cdF} \leq   H( y, \xi,\Theta)   + Cy(x^2+ y^\beta),
$$
where
\be\label{defH}
\begin{split}
H(y, \xi, \Theta) & =\Bigl\{p \bigl[q(1-\beta)-\psi(y)\bigr]- q(p-1) \xi\Bigr\} \Theta  
\\
&\quad  - q(q-1)\xi^2+2q \bigl[q(1-\beta)-\psi(y)\bigr]\xi \\
&\quad  -q(1-\beta) \bigl(q(1-\beta)+1\bigr)+2q(1-\beta) \psi(y)\,.
\end{split}
\ee

\smallskip
\noindent{\hskip 8mm}{\bf   Subcase 2.1: } $p \bigl[q(1-\beta)-\psi(y)\bigr]- q(p-1) \xi \le 0$.
\medskip

Using $\xi\ge 0$ and $0\le\psi(y)\le y_1$ in $Q$, we first observe that this implies $H\leq q\,f(\xi)$,
where
$$
f(\xi):=  -  (q-1)\xi^2  +2 q(1-\beta)\xi- (1-\beta) \bigl(q(1-\beta)+1\bigr)+2 (1-\beta)y_1. \\
$$
Computing  the reduced discriminant of this trinomial we notice that
\begin{align*}
\Delta
&=q^2(1-\beta)^2-(q-1)(1-\beta) \bigl(q(1-\beta)+1\bigr)+2 (q-1)(1-\beta)y_1\\
&= (1-\beta)\bigl[1-q\beta+2 (q-1)y_1\bigr]<0
\end{align*}
provided $q>\frac 1\beta= p-1$ 
and $y_1$  is sufficiently small. Therefore $f(\xi) <0$ for every $\xi \, \ge 0$ and we have
$$
p \bigl[q(1-\beta)-\psi(y)\bigr]- q(p-1) \xi \le 0 \quad \Longrightarrow \quad  H\leq q\,\max_{\xi\ge 0} f(\xi) <0\,.
$$
Hence ${\mathcal P}J<0$ provided $y_1$ is sufficiently small.

\medskip
\noindent{\hskip 8mm}{\bf   Subcase 2.2: } $p \bigl[q(1-\beta)-\psi(y)\bigr]- q(p-1) \xi > 0$.
\medskip

Here, we use 
\rife{u1d} to estimate 
$$
\Theta \leq [c_p(1-\beta)]^{p-1} (1+ Cy^{1+\beta})^{p-1} = \beta (1+ Cy^{1+\beta})^{p-1}\le  \beta + Cy^{1+\beta}.
$$
Since $p \bigl[q(1-\beta)-\psi(y)\bigr]- q(p-1) \xi >0$, it follows that
\begin{align*}
H  & \leq \left\{p \bigl[q(1-\beta)-\psi(y)\bigr]- q(p-1) \xi\right\} \beta  
\\
&\quad  - q(q-1)\xi^2+2q \bigl(q(1-\beta)-\psi(y)\bigr)\xi\\
&\quad   -q(1-\beta) \bigl(q(1-\beta)+1\bigr)+2q(1-\beta) \psi(y) + C y^{1+\beta}
\end{align*}
which yields, using $p\beta=1+\beta$, 
\begin{align*}
H  & \leq q \left\{   \xi^2-\xi+\beta(1-\beta) - q\bigl[\xi-(1-\beta)\bigr]^2-2\psi(y) \bigl[\xi-(1-\beta)\bigr]  \right\}
\\
&\quad  - (1+\beta)\psi(y) + C y^{1+\beta}.
\end{align*}
Now we observe that, for any $y\in [0,y_1]$, the function
$$
\phi(\xi):= \xi^2-\xi+\beta(1-\beta) - q[\xi-(1-\beta)]^2-2\psi(y) [\xi-(1-\beta)]
$$
is concave and we have
\be\label{ineqphi1}
\phi'(1-\beta)= 1-2\beta-2\psi(y) \leq 0
\ee
since $\beta\geq \frac12$ ($\iff p\leq 3$). Therefore, we have
\be\label{ineqphi2}
\phi(\xi)\leq \phi(1-\beta)= 0 \quad\hbox{ for all $\xi\geq 1-\beta$}
\ee
and we conclude that
$$
H  \leq  - (1+\beta)\psi(y) + C y^{1+\beta},
$$
hence
\begin{align*}
{y^2{\mathcal P}J\over cdF} & \leq   - (1+\beta)\psi(y) + C y^{1+\beta}  +Cy(x^2+ y^\beta)
\\
& 
\leq - (1+\beta)\frac{ y}{1+ y} + Cy(x^2+ y^\beta) \leq 0
\end{align*}
provided $x_1$ and $y_1$ are taken sufficiently small.
The proof of Lemma \ref{pjneg} is complete.
\qed

\medskip

\noindent \textit{Continuation of proof of Proposition \ref{propJ}.}
\smallskip

\textbf{Step 4. Initial and boundary conditions for $J$.} 
First observe that, for each $T'<T$, we have
\begin{equation}\label{verif regul J}
u\leq Cy \ \ \ {\rm in}\ D\times [0,T']
\end{equation}
 for some $C=C(T')>0$. Since $\gamma<q$, we have in particular
\begin{equation}\label{regul J}
J\in C(\overline D\times [0,T))\cap C^{2,1}(D\times(0,T)).
\end{equation}
Also in view of (\ref{verif regul J}) and $\gamma<q-1$, the coefficient $a(x,y,t)$ of the operator $\mathcal{P}$ (cf.~(\ref{def a b})) satisfies
\begin{equation}\label{bdd a b}
\hbox{$a$ is bounded in $D\times(0,T')$ for each $T'<T$.}
\end{equation}

Next, since $w=u_x$ satisfies
\begin{eqnarray}\label{eq:u x}
    w_t-\Delta w=p|\nabla u|^{p-2}\nabla u\cdot\nabla w\ \ \ {\rm
    in}\  Q_T,
\end{eqnarray}
  and is nonnegative nontrivial in $ (0,L)^2\times(0,T)$, 
  by the maximum principle and after a time shift, we may assume that
\begin{equation}\label{res:u x negative}
    u_x<0\ \ \ {\rm
    in}\   (0,L)^2\times(0,T).
\end{equation} 

Let now $x_1,y_1,\sigma$ be given by Lemma \ref{pjneg} and assume $\sigma<T/2$ without loss of generality.
By (\ref{extension}), (\ref{eq:u x}), (\ref{res:u x negative}), (\ref{sys
main 1}) and Hopf's Lemma, there exist constants $c_1,c_2>0$ such that
\begin{eqnarray*}
    &&u_x\leq -c_1y\ \ \ {\rm on}\ \{x_1\}\times(0,y_1)\times(T/4,T),\\ 
        &&u_x\leq -c_1x\ \ \ {\rm on}\
        (0,x_1)\times\{y_1\}\times(T/4,T),\\ 
    &&u\leq c_2y\ \ \ {\rm on}\ \{x_1\}\times(0,y_1)\times(0,T).
\end{eqnarray*}
With the above estimates, we check  the function $J$ on the lateral boundary: if $y=y_1$, we have,
 by~\rife{bound6},
\be\label{cond:J boundary rectangle 3}
J(x,y_1,t)\leq  \bigl[-c_1+ky_1^{-\gamma} (1+ y_1)(c_p y_1^{1-\beta}+ c y_1^2)^q \bigr]x\leq 0
 \ \ \ {\rm on}\ (0,x_1)\times\{y_1\}\times(T/4,T),
\ee
if $k$ is sufficient small. If $x=x_1$, we have
\begin{eqnarray}
J(x_1,y,t)
&\leq& -c_1y+kx_1c_2^qy^{q-\gamma} (1+ y) \notag \\
&\leq&  \bigl[-c_1+kx_1c_2^qy_1^{q-\gamma-1} (1+ y_1)\bigr]y\leq 0
        \ \ \ {\rm on}\ \{x_1\}\times(0,y_1)\times(T/4,T), \label{cond:J boundary rectangle 4}
\end{eqnarray}
 if $k$ is sufficiently small,  since $\gamma= q(1-\beta)<q-1$. Moreover, we clearly have 
\begin{eqnarray}
    &&J(x,0,t)=0\ \ \ {\rm on}\ (0,x_1)\times\{0\}\times(T/4,T),\label{cond:J boundary rectangle 1}  \\
    &&J(0,y,t)=0\ \ \ {\rm on}\ \{0\}\times(0,y_1)\times(T/4,T).\label{cond:J boundary rectangle 2}  
\end{eqnarray}
Finally, we recall that there exists $c_3>0$ such that
\begin{equation}\label{claim corner}
    u_x(x,y,  T-\sigma)\leq -c_3xy\ \ \ {\rm in}\ D.
\end{equation}
(This is a parabolic version of ``Serrin's corner Lemma''; see \cite[p.512]{LS}).
Now (\ref{claim corner}) implies
\begin{eqnarray}
    J(x,y,  T-\sigma) &\leq& -c_3xy+kxc_2^qy^{q-\gamma} (1+ y) \notag  \\
        &\leq& \bigl[-c_3+kc_2^qy_1^{q-\gamma-1} (1+ y_1)\bigr]xy\leq 0\ \ \ {\rm in}\ D
        \label{cond:J initial rectangle}
\end{eqnarray}
  if $k$ is sufficient small, since $\gamma<q-1$. 
Then (\ref{res:J negative}) follows from 
Lemma \ref{pjneg}, (\ref{cond:J boundary rectangle 3})-(\ref{cond:J boundary rectangle 2}), (\ref{cond:J
initial rectangle}) and the maximum principle.
Note that the use of the maximum principle is justified in view of (\ref{regul J})   and (\ref{bdd a b})
(or, alternatively, of the fact that $a<0$).

The proof of Proposition \ref{propJ} is complete.
\qed

\vskip1em

  By combining Proposition \ref{propJ} and Lemma \ref{barrier}, we shall now prove the upper estimate.

\smallskip
\noindent 
\textit{Proof of Theorem 1.1: the upper estimate in (\ref{mainEstimate}).}
It suffices to prove it for $x>0$ sufficiently small (the case $x<0$ will follow by considering $u(-x,y,t)$).
By Proposition \ref{propJ}, we know now that,  for some $k>0$, $q>p-1$ and $t_0\in (0,T)$, we have
$$
u_x+kxy^{-q(1-\beta)}(1+y)u^q\leq 0 \ \ \ {\rm in}\ Q=(0,x_1)\times(0,y_1)\times (  t_0,T).
$$
 Integrating over $(0,x)$ and using \rife{bound6}, we obtain 
\begin{equation}\label{ConclLemmeCore2}
\begin{split}
    u^{1-q}(x,y,t) & \geq (q-1) k\frac {x^2}2 y^{-q(1-\beta)}+ u^{1-q}(0,y,t)
    \\ & \geq (q-1) k\frac {x^2}2 y^{-q(1-\beta)} + [c_p y^{1-\beta}+ c y^2]^{1-q}  \ \ \ {\rm in}\ Q.
    \end{split}
\end{equation}
Starting with this estimate, we shall now apply the regularizing barrier lemma (Lemma \ref{barrier}).

Fix $x_0\in (0,x_1)$. Let $\eta\in(0,\eta_0)$, where $\eta_0$ is given by Lemma \ref{barrier}, and set
\begin{equation}\label{defrdD}
r=\displaystyle\frac {x_0}2, \quad d=\eta x_0,\quad D= (x_0-r, x_0+r)\times(0,d).
\end{equation}
Next we recall $\kappa$, given by Lemma \ref{barrier}:
\begin{equation}\label{defKappa}
\kappa=C_0(p) \,\eta^{1-\beta}(r^2+ T-t_0).
\end{equation}
We shall also use the notation $\tau=t-t_0$.

We claim now that there exists $\eta\in (0,\eta_0)$ such that, for any $x_0$ sufficiently small, we have
\begin{align}\label{compuz0}
\begin{cases}
&c_p y^{1-\beta} -\kappa  \frac{y^2}2>0, \\
\noalign{\vskip 2mm}
&\bigl[c_p y^{1-\beta} -\kappa  \frac{y^2}2\bigr]^{1-q}-\bigl[c_p y^{1-\beta}+ c y^2\bigr]^{1-q}
\leq (q-1) k\frac {x^2}2 y^{-q(1-\beta)}
\end{cases}
\qquad\hbox{  in $\overline D\times [t_0,T)$}
\end{align}
and 
\begin{align}\label{compuz}
\begin{cases}
&c_p \bigl[(d+ \eta \tau^{\frac1{1-\beta}} r^{\frac2{1-\beta}})^{1-\beta}-\eta^{1-\beta} \tau r^2 \bigr]
-\kappa  \frac{d^2}2>0, \\
\noalign{\vskip 2mm}
&\left\{c_p \bigl[(d+ \eta \tau^{\frac1{1-\beta}} r^{\frac2{1-\beta}})^{1-\beta}-\eta^{1-\beta}\tau r^2 \bigr]
-\kappa  \frac{d^2}2\right\}^{1-q}  \\
\noalign{\vskip 2mm}
&\qquad\qquad-\bigl[c_p d^{1-\beta}+ c d^2\bigr]^{1-q} \leq (q-1) k\frac {x^2}2 d^{-q(1-\beta)}
\end{cases}
\qquad\hbox{  in $[x_0-r, x_0+r]\times [t_0,T)$.} 
\end{align}
\vskip0.5em

Assume for the moment that \rife{compuz0}-\rife{compuz} hold; together with \rife{ConclLemmeCore2}, this implies
$$
u^{1-q}(x,y,t)    \geq \left\{c_p  y^{1-\beta} -\kappa  \frac{y^2}2\right\}^{1-q} 
 \quad\hbox{  in $\overline D\times  [t_0,T)$}
$$
and
$$
u^{1-q}(x,d,t)  \geq\left\{c_p [d+ \eta \tau^{\frac1{1-\beta}} r^{\frac2{1-\beta}})^{1-\beta}-\eta^{1-\beta} \tau r^2]-\kappa  \frac{d^2}2\right\}^{1-q}
 \quad\hbox{  in $[x_0-r, x_0+r]\times  [t_0,T)$,}
$$
and so both \rife{hypcontroluc0} and \rife{toit} will hold. 
We may then apply Lemma \ref{barrier} to deduce
\begin{equation}\label{conclusionUpper}
u_y(x_0,0,t) \le  d_p \eta^{-\beta} \bigl(  \tau r^2 \bigr)^{-\frac\beta{1-\beta}} \quad\hbox{  in $(t_0,T)$.}
\end{equation}
At time $t=T$ this gives
$$
u_y(x_0,0,T )\le C \eta^{-\beta}\, x_0^{-\frac{2\beta}{1-\beta}} = C \eta^{-\beta}\, x_0^{-\frac2{p-2}}
$$
which, jointly with (\ref{bound5a}), implies the upper estimate in (\ref{mainEstimate}).
\vskip0.5em
To conclude, we are thus left to prove \rife{compuz0}-\rife{compuz}. We note right away that the
first part of \rife{compuz0} is true whenever  $d$ is small enough, a  condition which holds as soon as $x_0$ is small enough (independently of $\eta$).
Similarly, observe that
if $x_0$ --~hence~$d$~-- is sufficiently small (independently of $\eta$), then
\begin{align*}
&\Bigl[c_p y^{1-\beta} -\kappa  \frac{y^2}2\Bigr]^{1-q}-\bigl[c_p y^{1-\beta}+ c y^2\bigr]^{1-q} \\
&\qquad\qquad  \leq  (c_p y^{1-\beta})^{1-q}\left[\bigl(1 - \frac\kappa2 c_p^{-1} y^{1+\beta}\bigr)^{1-q}-\bigl(1 +c c_p^{-1} y^{1+\beta}\bigr)^{1-q}\right] \\
&\qquad\qquad  \leq  C y^{(1-\beta)(1-q)}y^{1+\beta}=C y^{2-q(1-\beta)} \leq C d^2 y^{-q(1-\beta)}=C \eta^2 x_0^2\, y^{-q(1-\beta)}
\end{align*}
in $D\times (t_0,T)$. Here and in the rest of the proof, $C$ denotes a generic constant independent of $x_0$ and $\eta$.
Consequently, \rife{compuz0} holds as soon as $\eta$ is sufficiently small.

In order to prove \rife{compuz}, setting $\zeta=\eta \tau^{\frac1{1-\beta}} r^{\frac2{1-\beta}}$, we write
\begin{align*} 
c_p \bigl[(d+ \eta \tau^{\frac1{1-\beta}} r^{\frac2{1-\beta}})^{1-\beta}-\eta^{1-\beta}\tau r^2\bigr]-\kappa  \frac{d^2}2
&= c_p d^{1-\beta} \left[\left( 1+ \frac{\zeta}d  \right)^{1-\beta}-\left(\frac\zeta d\right)^{1-\beta}- \frac\kappa{2c_p} d^{1+\beta}\right] \\
&\ge c_p d^{1-\beta} \left[1-\left(\tau^{\frac1{1-\beta}}r^{\frac{1+\beta}{1-\beta}}\right)^{1-\beta}- \frac\kappa{2c_p} (\eta r)^{1+\beta}\right].
\end{align*}
By (\ref{defKappa}), (\ref{defrdD}), it follows that
\begin{align}  \label{compare1}
c_p \bigl[(d+ \eta \tau^{\frac1{1-\beta}} r^{\frac2{1-\beta}})^{1-\beta}-\eta^{1-\beta}\tau r^2\bigr]-\kappa \frac{d^2}2
\ge c_p d^{1-\beta} \left[1-Cx_0^{1+\beta}\left(T+(x_0^2+T)\eta_0^2 \right) \right]>0,
\end{align}
provided
\begin{align} \label{smallx0}
Cx_0^{1+\beta}\left(T+(x_0^2+T)\eta_0^2 \right) < 1.
\end{align}
Note that (\ref{smallx0}) is true for $x_0$ sufficiently small, which guarantees the first part of \rife{compuz}.

Next notice that
the convexity inequality $ (a+b)^{1-q}\geq a^{1-q}+ (1-q)a^{-q}b$ implies
$$
[c_p d^{1-\beta}+ c d^2]^{1-q}\geq (c_p d^{1-\beta})^{1-q} - C d^{-q(1-\beta)+2}.
$$
Combining this with (\ref{compare1}) and (\ref{smallx0}), it follows that
\begin{align*} &
\left\{c_p \bigl[(d+ \eta \tau^{\frac1{1-\beta}} r^{\frac2{1-\beta}})^{1-\beta}-\eta^{1-\beta}\tau r^2\bigr]-\kappa  \frac{d^2}2\right\}^{1-q}-\bigl[c_p d^{1-\beta}+ c d^2\bigr]^{1-q} \\
& \qquad \leq (c_p d^{1-\beta})^{1-q} \left\{ \left[1-Cx_0^{1+\beta}\left(T+(x_0^2+T)\eta_0^2 \right) \right]^{1-q} -1\right\} + C d^{-q(1-\beta)+2} \\
&  \qquad \leq  C d^{(1-\beta)(1-q)} x_0^{1+\beta}\left(T+(x_0^2+T)\eta_0^2 \right) + C d^{-q(1-\beta)+2} \\
&  \qquad \leq  C \eta^{1-\beta}\left(T+(x_0^2+T)\eta_0^2+ \eta_0^{1+\beta} \right)  x_0^2\,d^{-q(1- \beta)}. 
\end{align*}
Therefore, \rife{compuz} is satisfied as soon as $\eta$ is chosen sufficiently small. \qed

\section{A heuristic explanation of the singularity exponents through quasi-stationary approximation}

A possible heuristic explanation of the appearance of the number $2/(p-2)$ in the tangential singularity profile (\ref{mainEstimate}) can be obtained
using the idea of {\bf quasi-stationary approximation} along the family of 1D steady states.

\medskip

Recall the following family of 1D steady states, given by the translates of the reference solution 
$$ V(y)=c_p y^{1-\beta},$$
i.e.
$$V_a(y)=V(y+a)-V(a),  \quad  y\ge 0,  \ a\ge 0.$$
These special solutions verify
$$-V_a''={V_a'}^p,  \qquad V_a(0)=0,   \qquad V_a'(0)=d_p a^{-\beta}.$$
The idea is then to look for an approximate solution obtained by modulating in $a$, or moving on the manifold of steady-states $(V_a)_{a\ge 0}$.
More precisely, we set $U=u_{approx}$ given by
\begin{equation}\label{defApproxSol}
U(x,y,t) = V(y+h(t,x))-V(h(t,x)),
\end{equation}
which amounts to parametrize the solution by $a=h(t,x)$.
In particular, we have $U_y(x,y,t) = V'(y+h(t,x)).$
The function $h(t,x)$ is positive for $t<T$ and must satisfy $h(T,0)=0$ 
so that $U_y(0,0,T)=\infty.$

Note that this Ansatz means in some sense that
$-u_{yy} \sim (u_y)^p$ and $u_t \sim u_{xx}$ near the singularity,
already giving a rough clue to the parabolic nature of the scaling 
of the profiles in $t$ and $x$. 

With the above Ansatz, one has an interpretation of the lower estimate of the tangential profile in (\ref{mainEstimate})
as being {\bf a consequence of the constraint $U_{xx} \ge -C$, which comes from the maximum principle} (cf.~Proposition \ref{prop:maxprin}).
 Indeed, $U_{xx} \ge -C$ and $U_x(0,y,t)=0$ imply that
\begin{equation}\label{yrate0}
U_x \ge -Cx,\quad x>0.
\end{equation}
For $t<T$, restricting without loss of generality to $x>0$, we note that 
$$U_x(x,y,t)=d_p\bigl[(y+h(t,x))^{-\beta}-h^{-\beta}(t,x)\bigr]h_x\ge -d_ph_xh^{-\beta}(t,x),$$
where we used $h_x>0$ due to $U_x<0$, and that
$$U_x(x,h(t,x),t)=-ch_xh^{-\beta}(t,x).$$
Consequently, (\ref{yrate0}) is equivalent to $h_xh^{-\beta}(t,x) \le Cx$. By integration, it follows that
$h^{1-\beta}(t,x)\le h^{1-\beta}(t,0)+Cx^2$. Letting $t\to T$, we get $h(T,x)\le Cx^{2/(1-\beta)}$, which leads to
\begin{equation}\label{yrate}
 U_y(x,0,T)=V'(h(T,x))=d_p (h(T,x))^{-\beta}\ge Cx^{-2\beta/(1-\beta)}=Cx^{-2/(p-2)}.
\end{equation}
The fact that the upper estimate in (\ref{mainEstimate}) is exactly of this order means that
the constraint $U_{xx} \ge -C$ is satisfied in a minimal way by the parabolic flow.

\medskip
The same analysis can be done with the
time rate as well and actually enables one to recover also the exponent $1/(p-2)$ 
of the time rate\footnote{this observation 
doesn't seem to have been made in previous work on the time rate.}, namely
$$
\|\nabla u(t)\|_\infty\sim (T-t)^{-1/(p-2)}
$$
(the lower estimate is always true -- see \cite{CG96,GH08,QS07} -- whereas the upper estimate is only known for monotone increasing solutions in 1D;
see \cite{GH08} and cf. also \cite{QS07}).
This time the essential constraint is $|U_t| \le C$ (cf.~Proposition \ref{prop:maxprin}).
Indeed, one can easily see that
$|U_t| \le C$ is equivalent to $|(V(h))_t| \le C$.
Since $h(T,0)=0$, 
we thus have $V(h(t,0))\le C(T-t)$,
i.e. $h(t,0)\le C(T-t)^{1/(1-\beta)}$, or
\begin{equation}\label{timerate}
U_y(0,0,t) 
=V'(h(t,0))=d_p (h(t,0))^{-\beta} \ge C(T-t)^{-\beta/(1-\beta)}=C(T-t)^{-1/(p-2)}.
\end{equation}
Since the rates (\ref{yrate}) and (\ref{timerate}) violate the self-similar structure, or natural scaling, of the equation
(cf.~Remark (a) in Section 1.3),
so one can say that the maximum principle here wins against self-similarity.

 \begin{rem}\label{remp3}
  For $p>3$, the proof of the upper estimate in (\ref{mainEstimate})
 fails at the level of inequalities \rife{ineqphi1}--\rife{ineqphi2}. 
Actually, it can be seen along the lines of the proof of Lemma \ref{pjneg}
 that ${\mathcal P}J>0$ in some regions near the singularity
(more precisely, where $yu_yu^{-1}\sim (1-\beta)_+$ and $yu_y^{p-1}\sim\beta$).
However this might be technical, and it is presently open whether or not
the actual behavior of $u$ changes for $p>3$.

As for the above heuristic argument, although it does not a priori make a difference between the ranges
$2<p\le 3$ or $p>3$, it is not clear if such an argument can suggest more than a lower estimate of the profile.
Indeed, we stress that the heuristic argument gives a justification of the lower estimate in (\ref{mainEstimate})
only in view of the one-sided estimate 
\begin{equation}\label{onesided}
u_{xx}\ge -C.
\end{equation}
For $p\le 3$, our proof of the upper estimate also shows that the estimate $u_{xx}\ge -C$
is really optimal and that $u_{xx}$ is discontinuous near the singularity at $t=T$ (cf.~Remark~\ref{rem51}). 
If one could show for $p>3$ that $u_{xx}$ remains continuous (i.e., has a zero limit) 
near the singularity at $t=T$,
then the proof of Theorem \ref{thm:GBUprofileLower}  (as well as the heuristic argument),
 would imply that the final profile of $u_y$
is more singular than (\ref{mainEstimate}).

However, such property of $u_{xx}$ should be rather unstable and proving this might be quite delicate. 
Indeed, for any $p>2$ and {\it any} $\alpha\ge (p-1)/(p-2)$, 
a simple computation shows that the function
$$u(x,y,t)=c_p\bigl[(|x|^{2\alpha}+(T-t)^\alpha+y)^{1-\beta}-(|x|^{2\alpha}+(T-t)^\alpha)^{1-\beta}\bigr],$$
modeled after \rife{defApproxSol}, is a solution of  
 $$u_t-\Delta u=|\nabla u|^p+f$$ 
 in $Q=(-1,1)\times (0,1)\times (0,T)$
 with $u(x,0,t)=0$ and some $f\in L^\infty(Q)$. Moreover $u$ has an isolated
 gradient blowup point at $(0,0,T)$,
 $u$ satisfies \rife{onesided} and $|u_t|\le C$, but
 $u_{xx}$ is continuous near $(0,0,T)$ if and only if $\alpha> (p-1)/(p-2)$.
\end{rem}

\bigskip 
\noindent{\bf Acknowledgement.} \quad  
Part of this work was done during a visit of A.~Porretta at the Laboratoire Analyse, G\'{e}om\'{e}trie et Applications of Universit\'e Paris 13.
A.~Porretta wishes to thank these institutions for the invitation and the very kind hospitality given on this occasion.


\end{document}